\documentclass[a4paper,12pt,reqno]{amsart}

\usepackage{amssymb}

\numberwithin{equation}{section}

\theoremstyle{plain}
\newtheorem{thm}{Theorem}[section]
\newtheorem{cor}[thm]{Corollary}

\newtheorem{prop}[thm]{Proposition}
\newtheorem{conj}{Conjecture}

\theoremstyle{definition}
\newtheorem{defn}{Definition}
\newtheorem{example}{Example}
\newtheorem{rem}{Remark}

\def\C{\mathbb C}

\def\Z{\mathbb Z}

\def\vir{\operatorname{Vir}}

\def\Ann{\operatorname{Ann}}

\def\eps{\varepsilon}

\def\semiinfty{{\frac\infty2}}

\def\End{\operatorname{End}}

\DeclareMathOperator{\rank}{rank}

\def\Id{\operatorname{Id}}

\def\a{\alpha}

\def\sl{\mathfrak {sl}}
\def\slhat{{\widehat{\mathfrak {sl}}(2,\C)}}

\def\1{{\mathbf 1}}
\def\o{\otimes}

\def\g{\mathfrak g}
\def\ghat{\hat{\mathfrak g}}

\def\h{\mathfrak h}

\def\n{\mathfrak n}
\def\nhat{\hat{\mathfrak n}}

\def\ggbar{{\g \oplus \g}}
\def\ghatghat{{\ghat \oplus \ghat}}
\def\ghatKK{{\ghat_k \oplus \ghat_{\bar k}}}

\def\virCC{{\vir_c \oplus \vir_{\bar c}}}
\def\virghat{{\vir \ltimes \ghat}}
\def\virnil{{\vir \ltimes \hat{\mathfrak n}_+}}

\def\<{\langle}
\def\>{\rangle}
\def\({\left(}
\def\){\right)}

\def\pd#1{{\partial_{#1}}}
\def\equivt{\stackrel{\mod t^2}\equiv}

\begin{document}

\title[Regular representations, vertex algebras and semi-infinite cohomology]
{Modified regular representations\\ of affine and Virasoro algebras, VOA structure \\and semi-infinite cohomology.}

\author{Igor B. Frenkel}
\address{Department of Mathematics, Yale University\\ New Haven, CT 06520, USA}

\author{Konstantin Styrkas}
\address{Max-Planck-Institut f\"ur Mathematik\\ D-53111 Bonn, Germany}


\begin{abstract}
We find a counterpart of the classical fact that the regular representation $\mathfrak R(G)$ of a simple complex group $G$ 
is spanned by the matrix elements of all irreducible representations of $G$. 
Namely, the algebra of functions on the big cell $G_0 \subset G$ of the Bruhat decomposition
is spanned by matrix elements of big projective modules from the category $\mathcal O$ of representations
of the Lie algebra $\g$ of $G$, and has the structure of a $\ggbar$-module.

We extend both regular representations to the affine group $\hat G$, and show that 
the loop form of the Bruhat decomposition of $\hat G$ yields modified versions of $\mathfrak R(\hat G)$.
They involve pairings of positive and negative level modules, with the total value of the central charge required for the existence 
of non-trivial semi-infinite cohomology. In this paper we consider in detail the case $G=SL(2,\C)$, 
the corresponding finite-dimensional and affine Lie algebras, and the closely related to them Virasoro algebra.

Using the Fock space realization, we show that both types of modified regular representations for the affine and Virasoro algebras
become vertex operator algebras, whereas the ordinary regular representations have instead the structure of conformal field theories.
We identify the inherited algebra structure on the semi-infinite cohomology when the central charge is generic. We conjecture that 
for the integral values of the central charge the semi-infinite cohomology coincides with the Verlinde algebra and its counterpart 
associated with the big projective modules.
\end{abstract}

\maketitle

\setcounter{section}{-1}
\section{Introduction.}

The study of the regular representation of a simple complex Lie group $G$ is at the foundation of representation theory of $G$. Realized as the space of regular functions on $G$, the regular representation $\mathfrak R(G)$ carries two compatible structures of a $G$-bimodule and of a commutative associative algebra. An algebro-geometric version of the Peter-Weyl theorem establishes the decomposition of $\mathfrak R(G)$ into a direct sum of subspaces, spanned by matrix elements of all irreducible finite-dimensional representations $V_\lambda$ of $G$, indexed by integral dominant highest weights $\lambda \in \mathbf P^+$.
In other words, we have an isomorphism of $G$-bimodules
\begin{equation}\label{eq:Peter-Weyl classical}
\mathfrak R(G) = \bigoplus_{\lambda \in \mathbf P^+} V_\lambda \o V_\lambda^*.
\end{equation}
where $V_\lambda^*$ is the dual representation of $G$. The multiplication in $\mathfrak R(G)$ can be described in representation-theoretic terms as a pairing of intertwining operators for the left and right $\g$-actions with appropriate structural coefficients. Thus the algebra structure on $\mathfrak R(G)$ encodes the information about the tensor category of finite-dimensional $\g$-modules.

The representations of $G$ can also be viewed as modules over the simple complex Lie algebra $\g$ associated with $G$. In the case of the Lie algebra $\g$ it is natural to consider a larger collection of modules - namely, the Bernstein-Gelfand-Gelfand category $\mathcal O$. Infinite-dimensional $\g$-modules from the category $\mathcal O$ are not integrable, and therefore their matrix elements cannot be regarded as functions on $G$. However, one can interpret them as functions on the open dense subset $G^o \subset G$, given by the Gauss decomposition
\begin{equation}\label{eq:Gauss decomposition}
G^o = N_- \cdot T \cdot N_+,
\end{equation}
where $N_\pm$ is the upper and lower triangular unipotent subgroup of $G$, and $T$ is the diagonal maximal abelian subgroup.
The space $\mathfrak R(G^o)$ of regular functions on $G^o$ does not have the structure of a representation of $G$.
Nevertheless, the left and right infinitesimal actions of the Lie algebra $\g$ on this space are well-defined, 
and can be expressed in terms of explicit differential operators in the parameters of the Gauss decomposition \eqref{eq:Gauss decomposition}.
The enlarged regular representation $\mathfrak R(G^o)$ decomposes into the direct sum of bimodules spanned by the matrix coefficients of all "big" projective $\g$-modules $P_\lambda$, indexed by strictly anti-dominant highest weights $\lambda \in - \mathbf P^{++} = - (\mathbf P^+ + \rho)$, where $\rho$ is the half-sum of all positive roots of $\g$. Thus we obtain an isomorphism of $\g$-bimodules
\begin{equation}\label{eq:Peter-Weyl twisted}
\mathfrak R(G^o) \cong \bigoplus_{\lambda \in -\mathbf P^{++}} \( P_\lambda \o P_\lambda^* \) \biggr/ I_\lambda,
\end{equation}
where $P_\lambda^*$ is the module dual to $P_\lambda$, and $I_\lambda$ is the sub-bimodule of the matrix coefficients which vanish identically on the universal enveloping algebra $\mathcal U(\g)$.

It is important to note that the dual modules $P_\lambda^*$ do not belong to the category $\mathcal O$, but to its ``mirror image'', 
in which all highest weight modules are replaces by lowest weight modules. In order to stay in the category $\mathcal O$ 
we replace the open subset $G^o$ coming from the Gauss decomposition by the maximal cell in the Bruhat decomposition
\begin{equation}\label{eq:Bruhat decomposition}
G_0 = N_+ \cdot \mathbf{w_0} \cdot T \cdot N_+,
\end{equation}
where $\mathbf {w_0}$ is the longest element of the Weyl group $W$. Then we obtain a version of the isomorphism \eqref{eq:Peter-Weyl twisted},
\begin{equation}\label{eq:Peter-Weyl projective}
\mathfrak R(G_0) \cong \bigoplus_{\lambda \in -\mathbf P^{++}} \( P_\lambda \o P_\lambda^\star \) \biggr/ I_\lambda,
\end{equation}
where the 'twisted' duals $P_\lambda^\star$ differs from $P_\lambda^*$ by the automorphism $\omega$ of $\g$, which is induced by $\mathbf w_0$ and interchanges the positive and negative roots.

The theorems of Peter-Weyl type and the Gauss decomposition can be extended to the central extension of the loop group $\hat  G$ associated to $G$, and to the corresponding affine Lie algebra $\ghat$ and its universal enveloping algebra $\mathcal U(\ghat)$. In this infinite-dimensional case the space $\mathfrak R(\hat G)$ of regular functions on $\hat G$ is decomposed into the direct sum of subspaces $\mathfrak R_k(\hat G)$, corresponding to the value $k \in \mathbb Z$ of the central charge.
Using the version of the Gauss decomposition \eqref{eq:Gauss decomposition} known as the Birkhoff decomposition, one can show (see \cite{PS}) that for any $k \in \mathbb Z$ there is an isomorphism
\begin{equation}\label{eq:Peter-Weyl affine}
\mathfrak R_k(\hat G) \cong \bigoplus_{\lambda \in \mathbf P^k_+} \hat V_{\lambda,k} \o \hat V_{\lambda,k}^*,
\end{equation}
where $\lambda$ runs over the truncated alcove $\mathbf P^+_k \subset \mathbf P^+$, depending on $k$, and $\hat V_{\lambda,k}$ are the corresponding irreducible modules. Similarly, one can obtain decompositions of $\mathfrak R_k(\hat G^o)$ analogous to \eqref{eq:Peter-Weyl twisted}, where $\hat G^o$ is the maximal cell in the Birkhoff decomposition.
Viewing the decomposition \eqref{eq:Peter-Weyl affine} in terms of the Lie algebra $\ghat$ allows to extend it for all values of $k$, with
$\mathbf P^+_k = \mathbf P^+$ for $k \notin \mathbb Q$. 

To transform the dual module $\hat V_{\lambda,k}^*$ into a module from the category $\mathcal O$ for $\ghat$, one might apply again an automorphism of $\ghat$ which interchanges the positive and negative affine roots. However, it no longer belongs to the affine Weyl group, and the Bruhat decomposition for $\hat G$ does not have a maximal cell.
To overcome this obstacle we consider instead an intermediate between the Birkhoff and the affine Bruhat decompositions - the loop version of the finite-dimensional Bruhat decomposition, and the corresponding big cell
\begin{equation}
\label{eq:loop Bruhat decomposition} 
\hat G_0 = LN_+ \cdot \mathbf{w_0} \cdot \widehat{LT} \cdot LN_+,
\end{equation}
where $LN_\pm$ denote the loop groups with values in $\n_\pm$, and $\widehat{LT}$ is the central extension of the loop group with values in $T$.
The decomposition \eqref{eq:loop Bruhat decomposition} is especially useful for explicit realizations of the 
left and right regular $\ghat$-actions in terms of differential operators. However, we are still in ``semi-infinite'' distance 
from the category $\mathcal O$, and need to further apply a well-known procedure of ``changing the vacuum'', 
which has originated from the free field realizations of the Wakimoto modules and the irreducible representations $\hat V_{\lambda,k}$ 
(see \cite{FeFr,BF}). 
As a result of this procedure we obtain a modified affine version of the extended regular representation \eqref{eq:Peter-Weyl projective},
\begin{equation}\label{eq:Peter-Weyl affine twisted}
\mathfrak R'_k(\hat G_0) \cong \bigoplus_{\lambda \in -\mathbf P^{++}} 
\( \hat P_{\lambda,k-h^\vee} \o \hat P_{\lambda,-k-h^\vee}^\star \) \biggr/ \hat I_{\lambda,k},
\end{equation}
where $\hat P_{\lambda,k-h^\vee}$ and $\hat P^\star_{\lambda,-k-h^\vee}$ are the projective $\ghat$-modules and their 'twisted' duals, 
$\hat I_{\lambda,k}$ are appropriate sub-bimodules, and we assume $k \notin\mathbb Q$.
The levels are shifted by the dual Coxeter number $h^\vee$, so that the diagonal $\ghat$-action 
has the level $-2h^\vee$.

Like the Wakimoto modules, the bimodule $\mathfrak R'_k(\hat G_0)$ is realized as a certain Fock space, with two 
commuting $\ghat$-actions described explicitly. This realization is similar to the standard realization of the Wakimoto modules, 
but the actions of $\ghat$ contain a crucial new ingredient - the vertex operators, directly related to the screening operators 
used to construct intertwining operators for the affine Lie algebra.
We also establish that for $k \notin \mathbb Q$ the structure of the socle filtration of the non-semisimple bimodule $\mathfrak R'_k(\hat G_0)$ 
is the same as in the finite-dimensional case. In particular, $\mathfrak R'_k(\hat G_0)$ contains the distinguished sub-bimodule
\begin{equation}\label{eq:positive Peter-Weyl}
\mathfrak R'_k(\hat G) \cong \bigoplus_{\lambda \in \mathbf P^{+}} \hat V_{\lambda,k-h^\vee} \o \hat V_{\lambda,-k-h^\vee}^\star. 
\end{equation}
The shifts of the central charge by the dual Coxeter number no longer allow the interpretation of the bimodules in 
\eqref{eq:Peter-Weyl affine twisted} and \eqref{eq:positive Peter-Weyl} as spaces of matrix elements of $\ghat$-modules. 
Nevertheless, the structures of these bimodules are completely analogous to those of bimodules $\mathfrak R(G_0)$ and $\mathfrak R(G)$!

The vacuum module $\hat V_{0,k}$ of the affine Lie algebra $\ghat$ carries an extremely rich additional structure of a vertex operator algebra (VOA); 
other $\ghat$-modules $\hat V_{\lambda,k}$ become its representations (see \cite{FZ}). We show in this paper that the bimodules 
$\mathfrak R'_k(\hat G)$ and $\mathfrak R'_k(\hat G_0)$ also admit a vertex operator structure, compatible with $\ghat$-actions. 
In the proof we use the explicit Fock space realization of these bimodules. 
As in the finite-dimensional case, the VOA structure of the modified regular representations $\mathfrak R'_k(\hat G)$ and $\mathfrak R'_k(\hat G_0)$ 
encodes the information about fusion rules of the corresponding tensor categories of $\ghat$-modules.
Thus besides the vacuum modules there is a class of vertex operator algebras associated to affine Lie algebras with fixed central charges. 
It is also related to the algebras of chiral differential operators over the simple algebraic group $G$ recently studied in \cite{GMS} and \cite{AG}. 
A study of this relation might help to understand the geometric nature of  the modified regular representations.

On the other hand, the original regular representation $\mathfrak R_k(\hat G)$ does not seem to have a VOA structure. 
Instead it has the structure of a two-dimensional conformal field theory, which is an object of a different nature despite having local properties 
similar to those of a VOA. The bimodule $\mathfrak R_k(\hat G_0)$ has a structure of a generalized (non-semisimple!) conformal field theory.

It is well-known that the representation theory of $\ghat$ is closely related to the representation theory of the corresponding $\mathcal W$-algebra 
via the quantum Drinfeld-Sokolov reduction. In particular, one expects to have an analogue of the Peter-Weyl theorem for the $\mathcal W$-algebras. 
In this paper we consider in detail the simplest case of $\g = \sl(2,\C)$, when the corresponding $\mathcal W$-algebra 
is the infinite-dimensional Virasoro algebra. We give explicit realizations of the Virasoro bimodules, 
analogous to \eqref{eq:Peter-Weyl affine twisted} and \eqref{eq:positive Peter-Weyl}, and equip them with compatible VOA structures. 
The structures of the non-semisimple modified regular representations are quite parallel in all cases; we fully describe their socle filtrations.
Generalizations of our constructions to higher rank Lie algebras are straightforward, but their $\mathcal W$-algebra versions require more technicalities;
full details will be presented in a subsequent paper. 

A remarkable feature of all the modified bimodules that appear in the decompositions of Peter-Weyl type is that the central charge of the diagonal subalgebra 
is always equal to the special values that appear in the semi-infinite cohomology theory \cite{Fe,FGZ} - namely, $-2h^\vee$ for the affine Lie algebras 
and 26 for Virasoro. Moreover, thanks to a general result of \cite{LZ}, the corresponding semi-infinite cohomology spaces inherit a VOA structure 
from the modified regular representations. In our case, they degenerate into commutative associative superalgebras, and for generic central charge 
we establish isomorphisms between cohomology groups with coefficients in the corresponding modified regular representations of the affine and 
Virasoro algebras and their finite-dimensional counterparts. In particular, we show that the 0th semi-infinite cohomologies of the 
affine and Virasoro algebras are isomorphic to the Grothendieck ring of finite-dimensional representations of $G$. 
We conjecture that for integral $k$ they lead to the Verlinde algebra and its projective counterpart.

This paper is organized as follows.
In Section 1 we consider the Bruhat decomposition and the Peter-Weyl theorems in the finite-dimensional case with $G=SL(2,\C)$.
We give a Fock space realization of the algebra $\mathfrak R(G_0)$, and obtain explicit formulas 
for the $\g$-actions and decomposition theorems, which will later be used as prototypes of the infinite-dimensional case. 
In the last subsection we compute the Lie algebra cohomology with coefficients in $\mathfrak R(G)$ and $\mathfrak R(G_0)$.
In Section 2 we study the affine case, and use the loop version of the finite-dimensional Bruhat decomposition to obtain the 
modified Peter-Weyl theorems for the spaces $\mathfrak R'_k(\hat G)$ and $\mathfrak R'_k(\hat G_0)$. The Fock space realization of these spaces
equips them with VOA structures compatible with the regular $\ghat$-actions. The semi-infinite cohomology of $\ghat$ with coefficients 
in the modified regular representations $\mathfrak R'_k(\hat G)$ and  $\mathfrak R'_k(\hat G_0)$ for generic central charge is shown
to be isomorphic to its finite-dimensional counterpart.
In Section 3 we construct the analogues of the modified regular representations of the Virasoro algebra using the quantum Drinfeld-Sokolov
reduction. We compute the corresponding semi-infinite cohomology groups using methods developed in string theory, and prove that they are 
isomorphic to their affine counterparts.
Finally, in Section 4 we describe another class of vertex operator algebras 
obtained by the pairing of $\slhat$ and Virasoro modules. We also discuss generalizations of our results to Lie algebras of other types,
and to the integral values of the central charge. We conclude with conjectures on relations of the semi-infinite cohomology of $\mathfrak R'_k(\hat G)$ 
and $\mathfrak R'_k(\hat G_0)$ for $k \in \Z_{>0}$ with the Verlinde algebra, its projective counterpart and twisted equivariant K-theory.

We wish to thank G. Zuckerman for sharing his expertise on semi-infinite cohomology, and S. Arkhipov, F. Malikov for valuable comments. 
I.B.F. is supported in part by NSF grant DMS-0070551.

\section{Regular representation of $\sl(2,\C)$ on the big cell.}

\subsection{Regular representations of $\sl(2,\C)$.} 
Let $G = SL(2,\C)$. 
We define the left and right regular actions of $G$ on the space $\mathfrak R(G)$ of regular functions on $G$ by 
\begin{equation} \label{eq:group regular actions}
(\pi_l(g)\psi)(h) = \psi(g^{-1} h), \qquad
(\pi_r(g)\psi)(h) = \psi(h \, g), \qquad
g,h \in G.
\end{equation}
The multiplication in $\mathfrak R(G)$ intertwines both left and right regular actions.

Let $T, N_+$ denote the diagonal and unipotent upper-triangular subgroups of $G$.
The group $W = \operatorname{Norm}(T)/T$ is called the Weyl group. 
The Bruhat decomposition $G = N_+ \cdot W \cdot T \cdot N_+$ implies that every $g \in G$ can be factored as
$g = n \cdot w \cdot t \cdot n'$ for some $n,n' \in N_+, t \in T, w \in W$.

We denote by $G_0$ the big cell of the Bruhat decomposition, corresponding to 
the longest Weyl group element ${\mathbf w_0}$. Explicitly, $G_0$ is the dense open subset of $G$, consisting of $g \in G,$
\begin{equation}\label{eq:explicit Bruhat decomposition}
g = 
\begin{pmatrix} 1 && x \\ 0 && 1 \end{pmatrix}
\begin{pmatrix} 0 && -1 \\ 1 && 0 \end{pmatrix}
\begin{pmatrix} \zeta && 0 \\ 0 && \zeta^{-1} \end{pmatrix}
\begin{pmatrix} 1 && y \\ 0 && 1 \end{pmatrix}
\end{equation}
for some $x,y \in \C$ and $\zeta \in \C^\times$. The variables $x,y,\zeta$ can be viewed as coordinates on $G_0$, and thus
the algebra $\mathfrak R(G_0)$ of regular functions on $G_0$ is identified with the space $\C[x,y,\zeta^{\pm1}]$.

Let $\g = \sl(2,\C)$ be the Lie algebra of $G$, with the standard basis
\begin{equation*}
\mathbf e = \begin{pmatrix} 0 && 1 \\ 0 && 0 \end{pmatrix}, \qquad
\mathbf h = \begin{pmatrix} 1 && 0 \\ 0 && -1 \end{pmatrix}, \qquad
\mathbf f = \begin{pmatrix} 0 && 0 \\ 1 && 0 \end{pmatrix},
\end{equation*}
satisfying the commutation relations
$[\mathbf h, \mathbf e] = 2 \mathbf e, \ 
[\mathbf h, \mathbf f] = -2 \mathbf f, \ 
[\mathbf e, \mathbf f] = \mathbf h.$
The nilpotent subalgebras $\n_{\pm}$ and the Cartan subalgebra $\h$ of $\g$ are defined by
$\n_+ = \C \mathbf e, \, \h = \C \mathbf h, \, \n_- = \C \mathbf f.$
The element $\mathbf {w_0} \in W$ determines a Lie algebra involution $\omega$ of $\g$, such that $\omega(\n_\pm) = \n_\mp$ and $\omega(\h) = \h$,
defined by
\begin{equation}\label{eq:Cartan automorphism}
\omega(\mathbf e) = -\mathbf f, \qquad
\omega(\mathbf h) = -\mathbf h, \qquad
\omega(\mathbf f) = -\mathbf e.
\end{equation}

The infinitesimal regular actions of $\g$ on $\mathfrak R(G)$, corresponding to \eqref{eq:group regular actions}, are given by 
\begin{equation} \label{eq:algebra regular actions}
(\pi_l(x) \psi)(g) = \frac d{dt} \psi(e^{-t \, x} g) \biggr|_{t=0}, \qquad
(\pi_r(x) \psi)(g) = \frac d{dt} \psi(g \, e^{t \, x} ) \biggr|_{t=0}, \qquad x \in \g, \, g \in G.
\end{equation}

These formulas also define left and right infinitesimal actions of $\g$ on the space $\mathfrak R(G_0)$.
(These actions cannot be lifted to the group $G$,
because $G_0$ is not invariant under left and right shifts).
Elementary calculations yield the following explicit description of the regular $\g$-actions (cf. \cite{FP}).

\begin{prop} \label{thm:classical action}
The regular $\g$-actions on $\mathfrak R(G_0)$ are given by
\begin{equation}
\begin{split} \label {eq:classical left action}
\pi_l(\mathbf e) & =  -\pd x,\\
\pi_l(\mathbf h) & =  \zeta \pd \zeta - 2 x \pd x,\\
\pi_l(\mathbf f) & = - x \zeta \pd \zeta + x^2 \pd x + \zeta^{-2} \pd y.
\end{split}
\end{equation}
\begin{equation}
\begin{split} \label {eq:classical right action}
\pi_r(\mathbf e) & =  \pd y ,\\
\pi_r(\mathbf h) & =  \zeta \pd \zeta - 2 y \pd y,\\
\pi_r(\mathbf f) & =  y \zeta \pd \zeta - y^2 \pd y - \zeta^{-2} \pd x.
\end{split}
\end{equation}
\end{prop}

\subsection{Bosonic realizations}

We now reformulate the constructions of the previous section in terms of Fock modules for certain Heisenberg algebras.
These realizations admit generalizations to the affine and Virasoro cases, where the geometric approach to the regular representations becomes more subtle.

The operators $\beta = x, \ \gamma = -\pd x$ acting on polynomials in $y$ give a representation of the Heisenberg algebra with generators 
$\beta,\gamma$ and relation $[\beta,\gamma] = 1$. The polynomial space $\C[y]$ is then identified with its irreducible 
representation $F(\beta,\gamma)$, generated by a vector $\1$ satisfying $\gamma \, \1 = 0$. 

The operators $\bar\beta  = -y, \ \bar\gamma = \pd y$ generate a second Heisenberg algebra,
acting irreducibly in the space $F(\bar\beta,\bar\gamma) \cong \C[x]$. Here and everywhere else in this paper the 'bar' notation 
is used to denote the second copies of algebras and their generators; it does not denote the complex conjugation.

We identify $\h^* \cong \C$ so that $\mathbf P \cong \Z$. Whenever possible, we use the more
invariant notation in order to avoid possible numeric coincidences.
The operators $\1_\lambda =\zeta^\lambda$ and $a = \zeta \pd \zeta$ gives rise to the semi-direct product
$\C[a] \ltimes \C[\mathbf P]$, with $\C[a]$ acting on $\C[\mathbf P]$ by derivations: $a \1_\lambda = \lambda \1_\lambda$.

Thus, we get a realization of the algebra $\mathfrak R(G_0)$ of regular functions on $G_0$, with the $\ggbar$-action described by the abstract 
versions of the formulas \eqref{eq:classical left action}, \eqref{eq:classical right action}.

\begin{thm}\label{thm:classical bimodule action}
The space $\mathbb F= F(\beta,\gamma) \otimes F(\bar\beta,\bar\gamma) \otimes \C[\mathbf P]$ gives a realization of the algebra $\mathfrak R(G_0)$. In particular,
\begin{enumerate}
\item
The space $\mathbb F$ has a $\ggbar$-module structure, given by 
\begin{equation}\label{eq:classical boson left}
\begin{split}
\mathbf e & = \gamma, \\
\mathbf h & = 2 \, \beta \gamma + a, \\
\mathbf f & = -\beta^2 \gamma - \beta \, a +  \bar \gamma \, \1_{-2},
\end{split}
\end{equation}
\begin{equation}\label{eq:classical boson right}
\begin{split}
\bar {\mathbf e} & = \bar\gamma, \\
\bar {\mathbf h} & = 2 \, \bar\beta \bar\gamma + a, \\
\bar {\mathbf f} & = -\bar\beta^2 \bar\gamma - \bar\beta \, a +  \gamma \, \1_{-2}.
\end{split}
\end{equation}
\item
The space $\mathbb F$ has a compatible commutative algebra structure (i.e. the multiplication in $\mathbb F$ intertwines the $\ggbar$-action).
\end{enumerate}
\end{thm}

By specializing the action \eqref{eq:classical boson left} to the subspace $\ker \bar\gamma \subset \mathbb F$, we get the following well-known realizations of $\g$-action in the spaces $F_\lambda = F(\beta,\gamma) \o \C\1_\lambda$:
\begin{equation}
\begin{split} \label{eq:classical boson}
\mathbf e & = \gamma, \\
\mathbf h & = 2\beta \gamma + \lambda,\\
\mathbf f & = -\beta^2 \gamma - \lambda \, \beta.
\end{split}
\end{equation}

\begin{rem}
Simultaneous rescaling of the extra terms in \eqref{eq:classical boson left},\eqref{eq:classical boson right}, involving the shift $\1_{-2}$, by any multiple $\epsilon$ would preserve all the $\ggbar$ commutation relations.
For $\epsilon = 0$ such $\ggbar$-action degenerates into the product of two standard $\g$-actions \eqref{eq:classical boson}.
However, the multiplication in this
na\"{\i}ve bimodule loses much of its rich structure, and no longer encodes the information about the
fusion rules in the tensor category of finite-dimensional $\g$-modules.
\end{rem}

\subsection{$\ggbar$-module structure of the modified regular representation}
In this subsection we describe the socle filtration of the $\ggbar$-module $\mathbb F$.

For any $\lambda \in \h^*$, we denote by $V_\lambda$ the irreducible $\g$-module, generated by a highest weight vector $v_\lambda$ satisfying 
$\mathbf e \, v_\lambda = 0$ and $\mathbf h \, v_\lambda = \lambda \, v_\lambda$.

Recall that a $\g$-module $V$ is said to have a weight space decomposition, if
$$V = \bigoplus_{\mu\in\h^*} V[\mu], \qquad V[\mu] = \left\{v \in V \, \bigr| \, \mathbf h \, v = \mu \, v \right\}.$$
The restricted dual space $V' = \bigoplus_{\mu\in\h^*} V[\mu]'$ can be equipped with a $\g$-action, defined by
$$\<g \, v',  v\> = - \< v', \, \omega(g) \, v\>, \qquad g \in \g, \, v \in V, \, v' \in V',$$
where $\omega$ is as in \eqref{eq:Cartan automorphism}. We denote the resulting dual module $V^\star$.

We have an involution $\lambda \mapsto \lambda^\star$
of $\h^*$, determined by the condition $(V_\lambda)^\star \cong V_{\lambda^\star}$.
This involution can also be defined by $\lambda^\star = - \mathbf{w_0} (\lambda)$, where $\mathbf {w_0}$ is the
longest Weyl group element.

For $\g = \sl(2,\C)$, we have $\lambda^\star = \lambda$.
However, we keep the notation $\lambda^\star$, to indicate how our constructions 
generalize to Lie algebras of higher rank, where the involution is nontrivial.

\begin{thm}\label{thm:classical bimodule structure}
There exists a filtration 
\begin{equation}\label{eq:classical filtration}
0 \subset \mathbb F^{(0)} \subset \mathbb F^{(1)} \subset \mathbb F^{(2)} = \mathbb F
\end{equation}
of $\ggbar$-submodules of $\mathbb F$, such that
\begin{align}
\mathbb F^{(2)}/\mathbb F^{(1)} &\cong 
\bigoplus_{\lambda\in\mathbf P^+} V_{-\lambda-2} \o V^\star_{-\lambda-2} \label{eq:classical socle F2},\\
\mathbb F^{(1)}/\mathbb F^{(0)} &\cong 
\bigoplus_{\lambda\in\mathbf P^+} \( V_\lambda \o V^\star_{-\lambda-2} \oplus V_{-\lambda-2} \o V^\star_\lambda \) \label{eq:classical socle F1},\\
\mathbb F^{(0)} &\cong 
\bigoplus_{\lambda\in\mathbf P} V_\lambda \o V^\star_\lambda \label{eq:classical socle F0}.
\end{align}
\end{thm}

\begin{proof}
We introduce a filtration of $\ggbar$-submodules of $\mathbb F$
\begin{equation}
\label{eq:classical Fock filtration}
\dots \subset \mathbb F_{\le -2} \subset \mathbb F_{\le -1} \subset \mathbb F_{\le 0} \subset \mathbb F_{\le 1} \subset \mathbb F_{\le 2} \subset \dots,
\end{equation}
satisfying
$\bigcap_{\lambda\in\mathbf P} \ \mathbb F_{\le \lambda} = 0$ and 
$\bigcup_{\lambda\in\mathbf P} \ \mathbb F_{\le \lambda} = \mathbb F,$ where
$$\mathbb F_{\le \lambda} = F(\beta,\gamma) \o F(\bar\beta,\bar\gamma) \o \bigoplus_{\mu \le \lambda} \C\1_\mu, \qquad \lambda \in \mathbf P.$$
It is clear that $\mathbb F_{\le \lambda} / \mathbb F_{<\lambda} \cong F_\lambda \o F_{\lambda^\star}$; moreover, for $\lambda <0$
we have $F_\lambda \cong V_\lambda$, and for $\lambda \ge 0$ there is a short exact sequence
$0 \to V_\lambda\to F_\lambda \to V_{-\lambda-2}\to 0$. The linking principle for $\g$-modules 
implies that the successive quotients $\mathbb F_{\le \lambda} / \mathbb F_{<\lambda}$ and $\mathbb F_{\le \mu} / \mathbb F_{<\mu}$ of this filtration may be non-trivially linked only if $\mu = -\lambda-2$.

Thus we see that the $\ggbar$-module $\mathbb F$ splits into the direct sum of blocks
\begin{equation}\label{eq:classical double blocks}
\mathbb F = \mathbb F(-1) \oplus \bigoplus_{\lambda \in\mathbf P^+} \mathbb F(\lambda),
\end{equation}
where $\mathbb F(-1) \cong V_{-1} \o V^\star_{-1}$, and  $\mathbb F(\lambda)\cong 
\(V_{-\lambda-2} \o V^\star_{-\lambda-2}\) + \(F_\lambda \o F_{\lambda^\star}\)$
for $\lambda \in \mathbf P^+$; another way to obtain the decomposition \eqref{eq:classical double blocks} is by using the Casimir operator. 

It remains to describe the structure of $\mathbb F(\lambda)$ for each $\lambda \in \mathbf P^+$. By construction, $\mathbb F(\lambda)$ can be included in a short exact sequence
$0 \to V_{-\lambda-2} \o V^\star_{-\lambda-2} \to \mathbb F(\lambda) \to F_\lambda \o F_{\lambda^\star} \to 0.$
We conclude that there exists a filtration
$0 \subset \mathbb F(\lambda)^{(0)} \subset \mathbb F(\lambda)^{(1)} \subset \mathbb F(\lambda)^{(2)} = \mathbb F(\lambda),$
such that 
\begin{align*}\mathbb F(\lambda)^{(2)}/\mathbb F(\lambda)^{(1)} &\cong 
V_{-\lambda-2} \o V^\star_{-\lambda-2},\\
\mathbb F(\lambda)^{(1)}/\mathbb F(\lambda)^{(0)} &\cong 
\( V_\lambda \o V^\star_{-\lambda-2} \) \oplus \( V_{-\lambda-2} \o V^\star_\lambda \) ,\\
\mathbb F(\lambda)^{(0)} &\cong 
\( V_{-\lambda-2} \o V^\star_{-\lambda-2} \) + \( V_\lambda \o V^\star_\lambda \).
\end{align*}
In fact, the linking principle implies that the sum in $\mathbb F(\lambda)^{(0)}$ is direct:
$$\mathbb F(\lambda)^{(0)} \cong 
\( V_{-\lambda-2} \o V^\star_{-\lambda-2} \) \oplus \( V_\lambda \o V^\star_\lambda \).$$
Finally, we construct the filtration \eqref{eq:classical filtration} by setting
$$\mathbb F^{(0)} = \mathbb F(-1) \oplus \bigoplus_{\lambda \in \mathbf P^+} \mathbb F(\lambda)^{(0)},
\qquad
\mathbb F^{(1)} =  \mathbb F(-1) \oplus \bigoplus_{\lambda \in \mathbf P^+} \mathbb F(\lambda)^{(1)},$$
which obviously satisfies the required conditions 
\eqref{eq:classical socle F2},
\eqref{eq:classical socle F1},
\eqref{eq:classical socle F0}.
\end{proof}

\begin{rem}
For a Lie algebra $\g$ of higher rank, we will get a similar filtration of length $2\, l(\mathbf{w_0}) + 1$, and 
in addition to the regular blocks, corresponding to $\lambda \in \mathbf P^+$, and the most degenerate block $\mathbb F(-1)$,
there will be all intermediate types.
\end{rem}

The natural inclusion of algebras $\mathfrak R(G) \subset \mathfrak R(G_0)$ can be seen in the Fock space realizations.

\begin{cor}\label{thm:classical positive subalgebra}
There exists a subspace $\mathbf F \subset \mathbb F$ satisfying the following properties.
\begin{enumerate}
\item
$\mathbf F$ is a subalgebra of $\mathbb F$, and is generated by the elements from the submodule $V_1 \o V_1^\star$, corresponding to the matrix elements of the canonical representation of $G$.
\item
$\mathbf F$ is a $\ggbar$-submodule of $\mathbb F$, and is generated by the vectors $\{\1_\lambda\}_{\lambda \in \mathbf P^+}$. We have
\begin{equation}
\label{eq:classical Peter-Weyl}
\mathbf F = \bigoplus_{\lambda\in\mathbf P^+} \mathbf F(\lambda) \cong \bigoplus_{\lambda\in \mathbf P^+} V_\lambda \o V^\star_\lambda.
\end{equation}
\item
The space $\mathbf F$ is a realization of the algebra $\mathfrak R(G)$.

\end{enumerate}
\end{cor}

\noindent
In the polynomial realization, the generators of $\mathbf F$ from $V_1 \o V_1^\star$ are identified with functions
$$\psi_{11} = \zeta, \qquad \psi_{12} = x \zeta,  \qquad \psi_{21} = y \zeta, \qquad \psi_{22} = x y \zeta + \zeta^{-1},$$
which satisfy the relation $\psi_{11} \psi_{22} - \psi_{12} \psi_{21} = 1$. This establishes a very direct connection with the space of regular functions on the group $G = SL(2,\C)$.

\subsection {The generalized Peter-Weyl theorem}

In this section we interpret the space $\mathfrak R(G_0)$ of regular functions on $G_0$ and its Fock space realization $\mathbb F$ as the algebra of matrix elements of all modules from the category $\mathcal O$.

Recall that the Bernstein-Gelfand-Gelfand category $\mathcal O$ consists of all finitely generated, locally $\n_+$-nilpotent $\g$-modules. In particular, $V_\lambda\in\mathcal O$ for any $\lambda$.
If $V \in \mathcal O$, then $V^\star \in \mathcal O$.

For any $\g$-module $V$ we define $\mathbb M(V)$ to be the subspace of $\mathcal U(\g)'$, spanned by functionals 
\begin{equation}\label{eq:matrix element}
\phi_{v,v'}(x) = \<v', x \, v\>, \qquad v \in V, \, v' \in V', \, x \in \mathcal U(\g),
\end{equation}
where $\<\cdot,\cdot\>$ stands for the natural pairing between $V$ and $V'.$
The functionals \eqref{eq:matrix element} are called matrix elements of the representation $V.$

\begin{prop}\label{thm:matrix elements}
Introduce a $\ggbar$-module structure on the restricted dual $\mathcal U(\g)'$ by 
\begin{equation}\label{eq:dual enveloping bimodule}
(\pi_l(g) \phi) (x) = \phi(x g), \qquad\quad (\pi_r(g) \phi) (x) = - \phi(\omega(g) x)
\end{equation}
for any $\phi \in \mathcal U(\g)', \, g\in \g, \, x \in \mathcal U(\g)$. Then
\begin{enumerate}
\item
For any $\g$-module $V$, the space $\mathbb M(V)$ is a $\ggbar$-submodule of $\mathcal U(\g)'$.
\item
For any $\varphi \in \mathcal U(\g)',$ there exists a $\g$-module $V$, such that $\varphi \in \mathbb M(V).$ 
Moreover, if $\varphi$ is $\n_+ \oplus \n_+$-nilpotent, then $V$ can be chosen from the category $\mathcal O.$
\end{enumerate}
\end{prop}

\begin{proof}
To show that $\mathbb M(V)$ is invariant under the left action of $\g$, we compute
$$(\pi_l(g) \phi_{v,v'}) (x) = \phi_{v,v'}(x g) = \<v', x g \, v\> = \phi_{g v,v'}(x),$$
for any $x \in \mathcal U(\g),\, g \in \g,\, v \in V, \, v' \in V'.$ This shows that $y \phi_{v,v'}\in \mathbb M(V).$
The invariance under the right action follows from the computation
$$(\pi_r(g) \phi_{v,v'}) (x) = - \phi_{v,v'}(\omega(g) x) = - \<v', \omega(g) x \, v\> = \<y v', x \, v\> = \phi_{v, y v'}(x).$$

For the second part, assume $\varphi \in \mathcal U(\g)'$. Denote by $V$ the subspace of $\mathcal U(\g)'$, 
generated from $\varphi$ by the left action of $\g$.
Let $\varphi'$ be the restriction to $V$ of the unit $1 \in \mathcal U(\g) = \mathcal U(\g)''$. 
Equivalently, $\varphi'$ is the linear functional on $V$, determined by $\<\varphi', \psi\> = \psi(1),$ for any $\psi \in V \subset \mathcal U(\g)'$. 
We claim that $\varphi = \phi_{\varphi,\varphi'} \in \mathbb M(V)$. Indeed, for any $x \in \mathcal U(\g)$ we have
$$\phi_{\varphi,\varphi'} (x) = \<\varphi', x \varphi\> = (x\varphi)(1) = \varphi(x).$$

Finally, if $\varphi$ is left-$\n_+$-nilpotent, then $V$ is locally $\n_+$-nilpotent. 
Since $V$ is generated by a single element $\varphi,$ it belongs to category $\mathcal O.$ 
The right-$\n_+$-nilpotency condition guarantees that $\varphi'$ belongs to the {\it restricted} dual space $V'$.
\end{proof}

The elements of the universal enveloping algebra $\mathcal U(\g)$ may be regarded as the differential operators, acting on $\mathfrak R(G)$.
This gives an interpretation of the regular functions on $G$ (or even on $G_0$) as linear functionals on $\mathcal U(\g)$,
and thus to identifications of the spaces $\mathfrak R(G)$ and $\mathfrak R(G_0)$ with certain subspaces of $\mathcal U(\g)'$. In the explicit realizations
 $\mathbf F$ and $\mathbb F$ this correspondence is constructed using the algebraic analogue of the "co-unit" element of the Hopf algebra $\mathfrak R(G)$ - 
the linear functional $\<\cdot\>: \mathbb F \to \C$, defined by
\begin{equation}
\<\beta^m \bar\beta^n \1_\lambda\> = \delta_{m,0} \delta _{n,0}.
\end{equation}

\medskip

\begin{prop}
The linear map $\vartheta: \mathbb F \to \mathcal U(\g)'$, defined by $v \mapsto \vartheta_v$,
\begin{equation}\label{eq:map iota}
\vartheta_v(x) = \<\pi_l(x) v\>, \qquad v \in \mathbb F, \  x \in \mathcal U(\g).
\end{equation}
is an injective $\ggbar$-homomorphism.
\end{prop}

\begin{proof}
In terms of the polynomial realization, $\<\cdot\>$ corresponds to evaluating
a function $\psi(x,y,\zeta) \in \mathfrak R(G_0)$ at the element $\mathbf{w_0}$:
$\<\psi\> = \psi(0,0,1)$. This implies that for any $v \in \mathbb F$
\begin{equation}\label{eq:classical contravariant functional}
\<\mathbf e v\> = -\<\bar{\mathbf f} v\>, \qquad
\<\mathbf h v\> = -\<\bar{\mathbf h} v\>, \qquad
\<\mathbf f v\> = -\<\bar{\mathbf e} v\>.
\end{equation}
Therefore, for any $g \in \g$ and $x \in \mathcal U(\g)$ we have
$$\vartheta_{g v}(x) = \<x \, g v\> = \vartheta_v(x g) = (\pi_l(g)\vartheta_v)(x),$$
$$\vartheta_{\bar g v}(x) = \<x \, \bar g v\> = \<\bar g \, x v\> = - \<\omega(g) x v\> = 
- \vartheta_v(\omega(g) x) = (\pi_r(g)\vartheta_v)(x).$$
We conclude that the map $\vartheta$ is a $\ggbar$-homomorphism. To prove that it is injective, we need to show that for any nonzero $v \in \mathbb F$
there exists an element $x \in \ggbar$ such that $\<x v\> \ne 0.$

Since $\mathbb F$ is locally $\n_+$-nilpotent, we can pick $k\ge0$ such that
$\mathbf e^k v \ne 0,$ but $\mathbf e^{k+1} v = 0.$ Replacing $v$ by $\mathbf e^k v,$ we see that it suffices consider
the case of $v \ne 0$ such that $\mathbf e v = 0.$ Similarly, we may assume that $\bar{\mathbf e} v = 0.$

A vector $v$ satisfying $\mathbf e v = 0 = \bar{\mathbf e} v$ must have the form
$v = \sum_{\lambda\in\mathbf P} c_\lambda \1_\lambda$ with only finitely many $c_\lambda \ne 0.$
Using the formula for the Vandermonde determinant and the fact that
$$\<\mathbf h^m v\> = \sum_{\lambda\in\mathbf P} c_\lambda \lambda^m, \qquad m \ge 0,$$
we conclude that $\<\mathbf h^k v\> = 0$ for all $k\ge0$
if and only if all $c_\lambda$ vanish. Thus, $\theta_v = 0$ is equivalent to $v = 0$, which means that $\vartheta$ is an injection.
\end{proof}

The following statement is an algebraic version of the classical Peter-Weyl theorem.

\begin{thm}
The space $\mathfrak R(G)$ of regular functions on $G$ is spanned by the matrix elements of finite-dimensional irreducible $\g$-modules,
$$\mathfrak R(G) \cong \bigoplus_{\lambda\in \mathbf P^+} \mathbb M(V_\lambda).$$
The decomposition of $\mathfrak R(G)$ as a $\ggbar$-module is given by
\begin{equation*}
\mathfrak R(G) \cong \bigoplus_{\lambda\in \mathbf P^+} V_\lambda \o V^\star_\lambda.
\end{equation*}
\end{thm}

\begin{rem}
The subspace of $\mathcal U(\g)'$, corresponding to $\mathfrak R(G)$, is invariantly characterized as the restricted Hopf dual 
$\mathcal U(\g)'_{Hopf} \subset \mathcal U(\g)'$, defined by
$$\mathcal U(\g)'_{Hopf} = \{  \phi\in \mathcal U(\g)' \bigr| 
\exists \text{ two-sided ideal } J\subset \mathcal U(\g) \text { such that } \phi(J)=0 \text{ and }\operatorname{codim} J < \infty \}.$$
\end{rem}

The extended space $\mathfrak R(G_0)$ corresponds to a larger subalgebra of $\mathcal U(\g)'$, 
spanned by the matrix elements of all modules in the category $\mathcal O$.

Recall that the category $\mathcal O$ has enough projectives; we denote by $P_\lambda$ the indecomposable projective cover 
of the irreducible module $V_\lambda$. It is known that every indecomposable module in the category $\mathcal O$ with integral weights 
is isomorphic to a subfactor of the projective module, corresponding to some anti-dominant integral weight $\lambda$. 
In particular, this means that it suffices to consider the matrix elements of the big projective modules $\{P_\lambda\}_{\lambda<0}$.

The following result can be regarded as a non-semisimple generalization of the Peter-Weyl theorem.

\begin{thm}\label{thm:generalized Peter-Weyl classical}
The space $\mathfrak R(G_0)$ of regular functions on $G_0$ is spanned by the matrix elements of all big projective modules in the category $\mathcal O$,
$$\mathfrak R(G_0) \cong \bigoplus_{\lambda\in \mathbf P^+} \mathbb M(P_\lambda).$$
As a $\ggbar$-module, $\mathfrak R(G_0)$ is given by
$$\mathfrak R(G_0) \cong  \bigoplus_{\lambda\in-\mathbf P^{++}} \( P_{\lambda} \o P^\star_{\lambda} \) \biggr/ I_\lambda$$
where $I_\lambda$'s are the $\ggbar$-submodules of $P_{\lambda} \o P^\star_{\lambda}$, corresponding to identically vanishing matrix elements.
\end{thm}

\begin{proof}
We use the realization of $\mathfrak R(G_0)$ in the Fock space $\mathbb F$.
The inclusion \eqref{eq:map iota} provides the identification of $\mathbb F$ with a subspace of $\mathcal U(\g)'$.
Since $\mathbb F$ is locally $\n_+ \oplus \n_+$-nilpotent, Proposition \ref{thm:matrix elements} implies that for any
$v \in \mathbb F$ there exists a $\g$-module $W \in \mathcal O$ such that $\vartheta_v \in \mathbb M(W).$

Let $W = W_1 \oplus W_2 \oplus \dots \oplus W_m$ be the decomposition of $W$ into a direct sum of indecomposable
submodules. Each indecomposable component $W_i, \, i=1,\dots,m,$ is a subfactor of some big projective module $P_{\lambda_i}$.
Then $\mathbb M(W_i) \subset \mathbb M(P_{\lambda_i})$, and therefore we have
$$\mathbb M(W) = \mathbb M(W_1) + \mathbb M(W_2) + \dots + \mathbb M(W_m) \subset  \bigoplus_{\lambda\in -\mathbf P^{++}} \mathbb M(P_\lambda),$$
which shows that $\vartheta(\mathbb F) \subset  \bigoplus_{\lambda\in -\mathbf P^{++}} \mathbb M(P_\lambda).$ 
To prove that in fact $\vartheta(\mathbb F) = \bigoplus_{\lambda\in -\mathbf P^{++}} \mathbb M(P_\lambda)$, we compare the characters of the two spaces, and show that they have the same size.

For any $\lambda \in\mathbf P^+$ the $\ggbar$-module $\mathbb M(P_{-\lambda-2})$ is isomorphic to the quotient of the product $P_{-\lambda-2} \o P_{-\lambda-2}^\star$ by the kernel of the map 
\begin{equation}\label{eq:matrix elements map}
\Theta_\lambda: P_{-\lambda-2} \o P_{-\lambda-2}^\star \to \mathcal U(\g)', \qquad \Theta_\lambda(v \o v') = \phi_{v,v'}.
\end{equation}
Obviously, $I_\lambda = \ker \Theta_\lambda$ is a $\ggbar$-submodule of $P_{-\lambda-2} \o P_{-\lambda-2}^\star$; we describe it more explicitly.
It is known that the module $P_{-\lambda-2}$ has a filtration $0 \subset P^{(0)} \subset P^{(1)} \subset P_{-\lambda-2}$ such that
$$P^{(0)} \cong V_{-\lambda-2}, \qquad P^{(1)}/P^{(0)} \cong V_\lambda, \qquad P_{-\lambda-2}/P^{(1)} \cong V_{-\lambda-2},$$
and the dual filtration of the module $P_{-\lambda-2}^\star$ is given by
$$0 \subset \Ann(P^{(1)}) \subset \Ann(P^{(0)}) \subset P_{-\lambda-2}^\star.$$
They determine a filtration of the tensor product
\begin{equation*}
\begin{split}
0 \subset P^{(0)} \o  \Ann(P^{(1)}) \subset P^{(0)} \o  \Ann(P^{(0)}) + P^{(1)} \o  \Ann(P^{(1)}) \subset\\
\subset P^{(0)} \o P_{-\lambda-2}^\star + P^{(1)} \o  \Ann(P^{(0)}) + 
P_{-\lambda-2}\o  \Ann(P^{(1)}) \subset\\
\subset P^{(1)} \o  P_{-\lambda-2}^\star + P_{-\lambda-2} \o  \Ann(P^{(0)})
\subset P_{-\lambda-2} \o P_{-\lambda-2}^\star.
\end{split}
\end{equation*}
If $v\in P^{(0)}$ and $v' \in \Ann(P^{(0)})$, then $\phi_{v,v'}$ is the zero functional, since for any $x \in \mathcal U(\g)$ we have $x\, v \in P^{(0)}$ and $\phi_{v,v'}(x) = \<v', x \, v\> = 0$. Hence
the submodule $P^{(0)} \o  \Ann(P^{(0)})$ lies in the kernel of the map $\Theta_\lambda$, and similarly does
$P^{(1)} \o  \Ann(P^{(1)})$.
One can easily see that 
$$\Theta_\lambda \( P^{(1)} \o  \Ann(P^{(0)}) \) = \mathbb M(V_\lambda),$$
$$\Theta_\lambda \(  P^{(0)} \o P_{-\lambda-2}^\star \) = \Theta_\lambda \(  P_{-\lambda-2}\o  \Ann(P^{(1)})\) = \mathbb M(V_{-\lambda-2}).$$

It follows that the $\ggbar$-module $\mathbb M(P_{-\lambda-2})$ has a filtration
$$0 \subset \mathbb M^{(0)} \subset \mathbb M^{(1)} \subset \mathbb M^{(2)} = \mathbb M(P_{-\lambda-2})$$
such that
\begin{align*}
\mathbb M^{(2)}/\mathbb M^{(1)} & \cong  V_{-\lambda-2} \o V^\star_{-\lambda-2},\\
\mathbb M^{(1)}/\mathbb M^{(0)} & \cong  ( V_\lambda \o  V^\star_{-\lambda-2} ) \oplus  ( V_{-\lambda-2} \o  V^\star_\lambda ),\\
\mathbb M^{(0)} & \cong ( V_\lambda \o V^\star_\lambda ) \oplus ( V_{-\lambda-2} \o V^\star_{-\lambda-2} ).
\end{align*}

Thus, the block $\mathbb F(\lambda)$ of \eqref{eq:classical double blocks} is identified with the subspace, 
spanned by the matrix elements of the big projective module $P_{-\lambda-2}$.
Taking direct sums over all $\lambda \in \mathbf P^+$, adding the $\ggbar$-module $\mathbb M(P_{-1}) \cong V_{-1} \o V^\star_{-1}$, 
and comparing with Theorem \ref{thm:classical bimodule structure}, 
we see that $\bigoplus_{\lambda\in - \mathbf P^{++}} \mathbb M(P_\lambda)$ and $\mathbb F$ have the same characters. The statement of the theorem follows.
\end{proof}

\subsection{Cohomology of $\g$ with coefficients in regular representations}

The algebra $\mathfrak R(G)$ contains the subalgebra $\mathfrak R(G)^G$ of the conjugation-invariant functions on $G$,
which is linearly spanned by the characters of the irreducible finite-dimensional representations.
The subalgebra $\mathfrak R(G)^G$ is thus isomorphic to the Grothendieck ring of the finite-dimensional representations of $G$.

There is an isomorphism $\mathfrak R(G)^G \cong \C[\mathbf P]^W$, obtained by restricting the group characters to $\h$ and taking its Fourier expansion.
Finally, the algebra $\mathfrak R(G)^G$ also admits a cohomological interpretation, which will
be instrumental for further generalizations to the regular representations of the affine and Virasoro algebras. We briefly recall the definition of the cohomology of $\g$.

\begin{prop}
Let $\boldsymbol\Lambda = \bigwedge \g'$ be the exterior algebra of $\g'$ with unit $\1$. Then
\begin{enumerate}

\item The Clifford algebra, generated by $\{\iota(g), \, \eps(g')\}_{g \in \g, g' \in \g'}$ with relations
\begin{equation}
\label{eq:Clifford relations}
\{\iota(x), \iota(y) \} = \{ \eps(x'), \eps(y') \} = 0, \qquad
\{\iota(x), \eps(y')\} = \<y',x\>,
\end{equation}
acts irreducibly on $\boldsymbol\Lambda$, so that for any $\omega \in \boldsymbol\Lambda$ we have
$$\iota(g) \1 = 0, \qquad \eps(g') \omega = g'\wedge\omega, \qquad g \in \g,\, g' \in \g', \, \omega \in \boldsymbol \Lambda.$$

\item $\boldsymbol\Lambda$ is a commutative superalgebra,
$$\omega_1 \wedge \omega_2 = (-1)^{|\omega_1| \cdot |\omega_2|}\, \omega_2 \wedge \omega_1, \qquad \omega_1,\omega_2 \in \boldsymbol\Lambda,$$
where $|\cdot|$ is the natural grading on $\boldsymbol \Lambda$ satisfying $|\1| = 0, \  |\iota(g)| = -1, \ |\eps(g')| = 1.$

\item The $\g$-module structure on $\boldsymbol\Lambda$ is given by
$$\pi_{\boldsymbol\Lambda}(x) = \sum_i \eps(g'_i) \iota([g_i,x]),$$
where $\{g_i\}$ is any basis of $\g$, and $\{g'_j\}$ is the corresponding dual basis of $\g'$.

\end{enumerate}
\end{prop}

\begin{defn}
The cohomology $H^\bullet(\g;V)$ of $\g$ with coefficients in a $\g$-module $V$ is the cohomology of the graded complex
$C^\bullet(\g;V) = \boldsymbol\Lambda^\bullet \otimes V$, with the differential
\begin{equation}\label{eq:finite differential}
\mathbf d =\sum_i \eps(g'_i) \pi_V(g_i) - \frac12 \sum_{i,j} \eps(g'_i) \eps(g'_j) \iota([g_i,g_j]),
\end{equation}
where $\{g_i\}$ is any basis of $\g$, and $\{g'_i\}$ is the dual basis of $\g'$.
\end{defn}

The following is one of the fundamental results in Lie algebra cohomology, (see e.g. \cite{HS}).

\begin{thm}
\label{thm:vanishing classical cohomology}
For any finite-dimensional $\g$-module $V$ we have 
\begin{equation}
\label{eq:de Rham}
H^\bullet(\g; V) \cong V^\g \o H_{DR}^\bullet(G),
\end{equation}
where $H_{DR}^\bullet(G)$ denotes the holomorphic de Rham cohomology $H_{DR}^\bullet(G)$ of the Lie group $G$.
\end{thm}

If $V$ is a commutative algebra with a compatible $\g$-action, then its cohomology inherits the multiplication from $V$ and $\boldsymbol\Lambda$, 
and $H^\bullet(\g;V)$ becomes itself a commutative superalgebra. Moreover, the isomorphism \eqref{eq:de Rham} becomes an isomorphism of superalgebras,
with respect to the cup product in $H_{DR}^\bullet(G)$.

The diagonal $\g$-action in $\mathbf F$ corresponds to the coadjoint action of $G$ in $\mathfrak R(G)$; thus, we get

\begin{cor}
There is an isomorphism of commutative superalgebras
$$H^\bullet(\g;\mathbf F) = \C[\mathbf P]^W \o H_{DR}^\bullet(G).$$
\end{cor}

Our next goal is to study the cohomology of $\g$ with coefficients in the extended regular representation $\mathbb F \cong \mathfrak R(G_0)$. 
For infinite-dimensional $\g$-modules Theorem \ref{thm:vanishing classical cohomology} does not hold, and the cohomology $H^\bullet(\g;\mathbb F)$ does not reduce to $\mathbb F^\g \o H^\bullet_{DR}(G)$. We have instead

\begin{thm}
\label{thm:classical cohomology}
There is an isomorphism of commutative superalgebras 
$$H^\bullet(\g;\mathbb F) \cong \C[\mathbf P]^W \o \sideset{}{^\bullet}\bigwedge \C^2.$$
\end{thm}

\begin{proof}
It is easy to show using the results of \cite{W} that for $\lambda \ge -1$
$$H^n(\g;V_\lambda \o V^\star_{-\lambda-2}) = H^n(\g;V_{-\lambda-2} \o V^\star_\lambda) = 
\begin{cases}
\C, & n=1,2\\
0, & \text{otherwise}
\end{cases}$$
and that for $\lambda \ge0$ we have $H^n(\g;V_{-\lambda-2} \o V^\star_{-\lambda-2}) = 0 $ for all $n$.
The spectral sequence associated with the filtration of Theorem \ref{thm:classical bimodule structure}
can be used to show that
\begin{equation}
\label{eq:classical big cohomology}
H^n(\g;\mathbb F(-1)) =
\begin{cases}
\C, & n=1,2\\
0, & \text{otherwise}
\end{cases},\qquad
H^n(\g;\mathbb F(\lambda)) =
\begin{cases}
\C, & n=0,2\\
\C^2, & n=1\\
0, & \text{otherwise}
\end{cases}, \quad \lambda \ge0.
\end{equation}
Also, this spectral sequence shows that we have a natural isomorphism $H^0(\g;\mathbb F) \cong H^0(\g;\mathbf F).$

To explicitly get the generators of the commutative superalgebra $H^\bullet(\g;\mathbb F)$, we pick nonzero elements
$$\chi \in H^0(\g;\mathbb F(1)), \qquad 
\xi_{-1} \in H^1(\g;\mathbb F(-1)), \qquad
\eta_0 \in H^1(\g;\mathbb F(0)),$$
such that $\eta_0$ is not proportional to $\chi \, \xi_{-1}$. It is known that $H^0(\g;\mathbb F) \cong \C[\mathbf P]^W$ is isomorphic to the polynomial algebra $\C[\chi]$.
It is also clear that $H^\bullet(\g;\mathbb F)$ is a free $\C[\chi]$-module. For each $\lambda\ge0$, the set
$$B_{\le \lambda} = \{\xi_{-1}, \chi \xi_{-1}, \dots, \chi^{\lambda+1} \, \xi_{-1}\} \bigcup \{\eta_0,\chi \, \eta_0, \dots, \chi^\lambda \, \eta_0\}$$
consists of $2\lambda+3$ linearly independent elements, and in view of \eqref{eq:classical big cohomology} is a basis of $H^1(\g;\mathbb F_{\le \lambda})$.
Finally, one can check that $\eta_0 \, \xi_{-1} \ne 0$, and thus the elements $\{\eta_0 \, \xi_{-1}, \chi \, \eta_0 \, \xi_{-1},\dots, \chi^{\lambda+1}\, \eta_0 \, \xi_{-1}\}$
give a basis of $H^2(\g;\mathbb F_{\le\lambda})$ for each $\lambda \ge -1$.

It follows that $H^\bullet(\g;\mathbb F) \cong \C[\chi] \otimes \bigwedge^\bullet[\xi_{-1},\eta_0]$, and the theorem is proven.
\end{proof}

\begin{rem}
One of the ingredients in the exterior algebra part of the cohomology $H^\bullet(\g;\mathbb F)$ is the exterior algebra $\bigwedge^\bullet \h$, 
corresponding to $\bigwedge^\bullet[\eta_0]$ above. It would
be interesting to obtain an invariant characterization of the remaining part of $H^\bullet(\g;\mathbb F)$ for arbitrary $\g$.
\end{rem}

\begin{rem}
In each of the two-dimensional spaces $H^1(\g;\mathbb F(\lambda))$ there is a unique up to proportionality cohomology class $\xi_\lambda$ divisible by $\xi_{-1}$;
the elements $\frac {\xi_\lambda}{\xi_{-1}}$ constitute a basis of $H^0(\g;\mathbb F) \cong \C[\mathbf P]^W$, associated with
the characters of big projective modules (cf. \cite{La}).
\end{rem}

\section{Modified regular representations of the affine Lie algebra $\slhat$.}

\subsection{Regular representations of $\slhat$}
Let $\hat G$ be the central extension of the loop group $LG$, associated with $G=SL(2,\C)$ (see \cite{PS}), and let $\ghat$ be the corresponding Lie algebra.
As we discussed in the introduction, there is no maximal cell in the affine Bruhat decomposition, 
and thus we will use the loop version \eqref{eq:loop Bruhat decomposition} of the finite-dimensional one.
An additional advantage is that we get an explicit realization of the left and right regular $\ghat$-actions, 
analogous to the finite-dimensional case.

The standard basis of $\ghat$ consists of the elements $\{\mathbf e_n, \mathbf h_n, \mathbf f_n\}_{n \in \Z}$ and the central element $\mathbf k$, subject to the commutation relations
\begin{gather*}
[\mathbf h_m, \mathbf e_n] = 2 \mathbf e_{m+n}, \qquad
[\mathbf h_m, \mathbf f_n] = -2 \mathbf f_{m+n}, \qquad
[\mathbf h_m, \mathbf h_n] = 2 m \delta_{m+n,0} \mathbf k,\\
[\mathbf e_m, \mathbf f_n] = \mathbf h_{m+n} + m \, \delta_{m+n,0} \mathbf k, \qquad
[\mathbf e_m, \mathbf e_n] = [\mathbf f_m, \mathbf f_n] = 0.
\end{gather*}

The Lie algebra $\ghat$ has a $\Z$-grading $\ghat = \bigoplus_{n \in \Z} \ghat[n]$, determined by
$$\deg \mathbf f_n = \deg \mathbf h_n = \deg \mathbf e_n = -n, \qquad \deg \mathbf k = 0,$$
We introduce subalgebras $\ghat_\pm = \bigoplus_{\pm n>0} \g[n]$;
the finite-dimensional Lie algebra $\g$ is naturally identified with a subalgebra in $\ghat[0]$.

The element $\mathbf {w_0}$ of the classical Weyl group defines an involution $\hat\omega$ of $\ghat$, such that
\begin{equation}
\label{eq:affine Cartan automorphism}
\hat\omega(\mathbf e_n) = - \mathbf f_n, \quad
\hat\omega(\mathbf h_n) = - \mathbf h_n, \quad
\hat\omega(\mathbf f_n) = - \mathbf e_n, \quad
\hat\omega(\mathbf k) = - \mathbf k.
\end{equation}

We use the loop version of the finite-dimensional Bruhat decomposition \eqref{eq:Bruhat decomposition}, and factorize the central extension $\widehat{LT}$
into the product of loops that extend holomorphically inside and outside of the unit circle.
The analogue of \eqref{eq:explicit Bruhat decomposition} is the formal decomposition
$$g = \exp \( \sum_{n\in\Z} x_n \mathbf e_n \)  \  \mathbf{w_0}\,  \tau^\mathbf k \ 
\exp \( \sum_{m < 0} \zeta_m \mathbf h_m \) \zeta^{\mathbf h_0} \exp \( \sum_{m > 0} \zeta_m \mathbf h_m \) \exp \( \sum_{n\in\Z} y_n \mathbf e_n \).$$
The polynomial algebra $\mathfrak R_0(\hat G_0) = \C[\{x_n\}, \{y_n\}, \{\zeta_{n\ne0}\},\zeta^{\pm1}]$
can be thought of as the algebra of regular functions on the big cell of the loop group, and
for $\mathfrak R(\hat G_0)$ we get
$$\mathfrak R(\hat G_0) = \mathfrak R_0(\hat G_0) \o \C[\tau^{\pm1}] = \bigoplus_{\varkappa \in \Z} \mathfrak R_\varkappa(\hat G_0), \qquad
\mathfrak R_\varkappa(\hat G_0) = \mathfrak R_0(\hat G_0) \o \C\tau^\varkappa$$

Note that for each $\varkappa$ the subspace $\mathfrak R_\varkappa(\hat G_0)$ is a $\ghatghat$-submodule of $\mathfrak R(\hat G_0)$, but it is not a subalgebra of $\mathfrak R(\hat G_0)$ when $\varkappa \ne 0$ ! 
It is easy to see that the infinitesimal regular $\ghat$-actions
of the central element $\mathbf k$ on $\mathfrak R_\varkappa(\hat G_0)$ are given by
\begin{equation}
\label{eq:unmodified central charges}
\pi_l(\mathbf k) = - \varkappa \cdot \Id, \qquad\qquad \pi_r(\mathbf k) = \varkappa \cdot \Id.
\end{equation}
As vector spaces, all $\mathfrak R_\varkappa(\hat G_0)$ are identified with the same polynomial space,
and one can compute the infinitesimal regular actions of $\ghat$ by treating $\varkappa$ as a complex parameter. In particular, the
regular actions of $\ghat$ make sense for arbitrary $\varkappa \in \C$. Computations yield the following
description, analogous to Proposition \ref{thm:classical action}.

\begin{thm} \label{thm:affine action}
The regular actions of $\ghat$ on $\mathfrak R_\varkappa(\hat G_0)$ are given by \eqref{eq:unmodified central charges} and
\begin{equation} \label{eq:left affine action}
\begin{split}
\pi_l(\mathbf e_n) &= -\pd{x_n},\\
\pi_l(\mathbf h_n) &= - 2\sum_{i\in\Z} x_i \pd{i_{n+n}} +
\begin{cases}
\pd{\zeta_n} + 2n \varkappa\, \zeta_{-n}, & n>0\\
\zeta \, \pd {\zeta}, & n=0\\
\pd{\zeta_n}, & n<0 
\end{cases} ,\\
\pi_l(\mathbf f_n) &= \sum_{i,i' \in \Z} x_i x_{i'} \pd{x_{i+i'+n}} - 
\sum_{j<0} x_{j-n}\pd{\zeta_j} - x_{-n} \zeta\, \pd {\zeta} - \sum_{j>0} x_{j-n} \(\pd{\zeta_j} + 2 j \varkappa\, \zeta_{-j} \)  - \\
&- \varkappa \; n x_{-n} + \zeta^{-2} \sum_{j,j'>0} \mathrm s_{j'}(-2\zeta_1,-2\zeta_2,\dots) \, \mathrm s_j(-2\zeta_{-1},-2\zeta_{-2},\dots)  \pd{y_{n-j+j'}},
\end{split}
\end{equation}

\begin{equation} \label{eq:right affine action}
\begin{split}
\pi_r(\mathbf e_n) &= \pd{y_n},\\
\pi_r(\mathbf h_n) &= -2\sum_{i\in\Z} y_i \pd{y_{i+n}} +
\begin{cases} \pd{\zeta_n} & n>0, \\ \zeta \, \pd {\zeta}, & n=0,\\ \pd{\zeta_n} - 2 n \varkappa\, \zeta_{-n}, & n<0. \end{cases} ,\\
\pi_r(\mathbf f_n) &= - \sum_{i,i' \in \Z} y_i y_{i'} \pd{y_{i+i'+n}} + 
\sum_{j>0} y_{j-n}\pd{\zeta_j} + y_{-n} \zeta\, \pd {\zeta} + \sum_{j<0} y_{j-n} \(\pd{\zeta_j} - 2 j \varkappa\, \zeta_{-j} \)  - \\
&- \varkappa \, n y_{-n} - \zeta^{-2} \sum_{j,j'>0} \mathrm s_{j'}(-2\zeta_{-1},-2\zeta_{-2},\dots) \, \mathrm s_j(-2\zeta_1,-2\zeta_2,\dots)  \pd{x_{n+j-j'}},\\
\end{split}
\end{equation}
where the Schur polynomials $\mathrm s_k(\a_1,\a_2,\dots)$ are defined by
\begin{equation*}
\mathrm s_m(\a_1,\a_2,\dots) = \sum_{\substack {l_1,l_2, \ldots \ge 0\\ l_1+2l_2 + \ldots = m}} 
\frac {\a_1^{l_1} \a_2^{l_2} \dots}{l_1! l_2!\dots}.
\end{equation*}
\end{thm}

\begin{proof}
The presence of the central extension requires the use of some elementary cases of the Campbell-Hausdorff formula in our
computations;  we use the identity
$$\exp(B) \exp(tA) \equivt  \exp \(t \sum_{j=0}^\infty \frac 1{j!} \underbrace{[B,\dots,[B,[B,A]]\dots]}_{j \text{ commutators}}\) \exp(B).$$
For example, to derive the last of \eqref{eq:right affine action}, we use the formulas:
\begin{equation*}
\begin{split}
\exp \( \sum_{i\in\Z} y_i \mathbf e_i \)  \exp \( t \mathbf f_n \) \equivt &
\exp \(t \mathbf f_n \) 
\exp \(- t n y_{-n}\mathbf k + t \sum_{i\in\Z} y_i \mathbf h_{i+n}  \) \times \\
&\times \exp \(-t \sum_{i,i'\in\Z} y_i y_{i'} \mathbf e_{i+i'+n} \)
 \exp \( \sum_{i\in\Z} y_i \mathbf e_i \),\\
\exp \(\sum_{m>0} \zeta_m \mathbf h_m \) \exp \( t \mathbf f_n \) \equivt &
\exp \( t \sum_{j>0} \mathrm s_j(-2\zeta_1,-2\zeta_2,\dots) \mathbf f_{n+j} \)
\exp \( \sum_{m>0} \zeta_m \mathbf h_m \),\\
\zeta^{\mathbf h_0}  \exp \( t \mathbf f_n \) \equivt &
\exp \( t \zeta^{-2}\, \mathbf f_n \)  \zeta^{\mathbf h_0},\\
\exp \(\sum_{m<0} \zeta_m \mathbf h_m \) \exp \(t \mathbf f_n \) \equivt &
\exp \(t \sum_{j'>0} \mathrm s_{j'}(-2\zeta_{-1},-2\zeta_{-2},\dots) \mathbf f_{n-j'} \)
\exp \( \sum_{m<0} \zeta_m \mathbf h_m \),\\
\mathbf {w_0}  \exp \( t \mathbf f_n \) \equivt &
\exp \(- t \, \mathbf e_n \) \mathbf {w_0}.
\end{split}
\end{equation*}
Combining these equations, we get the desired formulas. We leave the technical calculations to the reader.
\end{proof}

\subsection{Vertex operator algebras: review and useful examples}

We aim to endow $\mathfrak R_\varkappa(\hat G_0)$ (or its modification) with a structure similar to that of an associative commutative algebra.
The relevant formalism is provided by the vertex algebra theory.

We recall the definitions of vertex and vertex operator algebras in the most convenient to us form. For more details and equivalent alternative definitions, we refer the reader to the books on the subject \cite{FLM,BFr}.

Let $\mathcal V$ be a vector space, equipped with a linear correspondence 
\begin{equation}\label{eq:state field correspondence}
v \mapsto \mathcal Y(v,z) = \sum_{n \in \Z} v_{(n)} z^{-n-1}, \qquad v_{(n)} \in \End(\mathcal V).
\end{equation}
We refer to such formal $\End(\mathcal V)$-valued generating functions as 'quantum fields'.

We say that $\mathcal V$ satisfies the locality property, if for any $a,b \in \mathcal V$
\begin{equation}\label{eq:locality}
(z-w)^N [\mathcal Y(a,z),\mathcal Y(b,w)] = 0  \quad \text{ for } N \gg 0 
\end{equation}
in the ring of $\End(\mathcal V)$-valued formal Laurent series in two variables $z,w$.

A vector $\1 \in \mathcal V$ is called the vacuum vector, if it satisfies
\begin{equation}\label{eq:vacuum}
\mathcal Y(\1,z) = \Id_{\mathcal V}, \qquad \mathcal Y(v,z)\1 \bigr|_{z=0} = v.
\end{equation}

An element $\mathcal D \in \End(\mathcal V)$, is called the infinitesimal translation operator, if it satisfies
\begin{equation}\label{eq:infinitesimal translation}
\qquad \mathcal D \, \1 = 0, \qquad\qquad [\mathcal D,\mathcal Y(v,z)] = \frac d{dz} \mathcal Y(v,z), \qquad \text{ for all } \   v \in \mathcal V.
\end{equation}

\begin{defn}
The space $\mathcal V$ is called a vertex algebra, if it is equipped with a linear map \eqref{eq:state field correspondence}, vacuum vector $\1$, and infinitesimal translation operator $\mathcal D$, 
satisfying the axioms \eqref{eq:locality}, \eqref{eq:vacuum}, \eqref{eq:infinitesimal translation} above. 
\end{defn}

Vertex superalgebras are defined as usual by inserting $\pm$ signs according to parity.
A vertex superalgebra $\mathcal V$ is called bi-graded, if it has $\Z$-gradings, $|\cdot|$ and $\deg$,
$$\mathcal V = \bigoplus_{m,n\in\Z} \mathcal V^m[n], \qquad 
\mathcal V^m[n] = \left\{ v \in\mathcal V \, \biggr| \, |v| = m \text{ and } \deg v = n \right\},$$
such that the parity in superalgebra is determined by $|\cdot|$, and for any homogeneous $v$
$$v \mapsto \mathcal Y(v,z) = \sum_{n\in\Z} v_{(n)} z^{-n-1}, \qquad \text { with } |v_{(n)}| = |v| \text{ and } \deg v_{(n)} = \deg v - n - 1 .$$
In particular, for the vacuum we must have \  $| \1 | = \deg \1 = 0$. Also, we write $|\mathcal Y(v,z)| = |v|$ and $\deg \mathcal Y(v,z) = \deg v$
for the quantum field $\mathcal Y(v,z)$, if the above conditions are satisfied.

A vertex algebra $\mathcal V$ is called a vertex operator algebra (VOA) of rank $c \in \C$, if there exists an element $\boldsymbol\omega \in \mathcal V$, usually called the Virasoro element, such that the operators $\{\mathcal L_n\}_{n \in \Z}$ defined by
$$\mathcal Y(\boldsymbol\omega,z) = \sum_{n \in \Z} \mathcal L_n z^{-n-2},$$
satisfy $\mathcal L_{-1} = \mathcal D$, and the Virasoro commutation relations
$$[\mathcal L_m,\mathcal L_n] = (m-n) \mathcal L_{m+n} + \delta_{m+n,0} \frac{m^3-m}{12}\, c.$$

We define the the normal ordered product of two quantum fields $X(z)$ and $Y(z)$ by
\begin{equation*}
:X(z)Y(w): = X_-(z)Y(w) + Y(w) X_+(z),
\end{equation*}
where $X_\pm(z)$ are the regular and principal parts of $X(z) = \sum_{n \in \Z} X_{(n)} z^{-n-1}$,
$$X_+(z) = \sum_{n\ge0} X_{(n)} z^{-n-1}, \qquad X_-(z) = \sum_{n < 0} X_{(n)} z^{-n-1}.$$
For products of three or more quantum fields, the normal ordered product is defined inductively, starting from the left.
In general, the normal ordered product is neither commutative nor associative.

The following 'reconstruction theorem' is an effective tool for constructing vertex algebras.

\begin{prop}\label{thm:free field construction}
Let $\mathcal V$ be a vector space with a distinguished vector $\1$ and a family of pairwise local $\End(\mathcal V)$-valued quantum fields $\{X^\a(z)= \sum_{n \in \Z}X^\a_{(n)} z^{-n-1}\}_{\a \in \mathfrak I}.$ 
Suppose $\mathcal V$ is generated from $\1$ by the action of the Laurent coefficients of quantum fields $X^\a(w)$, and that
the vectors $\{X^\a(z) \1 \bigr|_{z=0}\}_{\a \in \mathfrak I}$ are linearly independent in $\mathcal V$.
Then the operators
$$\mathcal Y\(X^{\a_1}_{(-n_1-1)}\dots X^{\a_k}_{(-n_k-1)} \1,z\) = 
\ :X^{\a_1}(z)^{(n_1)} \dots X^{\a_k}(z)^{(n_k)}:,$$
where $X(z)^{(n)} = \frac 1{n!} \frac {d^n}{dz^n} X(z)$,
satisfy \eqref{eq:locality} and  \eqref{eq:vacuum}.

If a linear operator $\mathcal D \in \End(\mathcal V)$ satisfies  $\mathcal D \1 = 0$ and 
$[\mathcal D,X^\a(z)] = \frac d{dz} X^\a(z)$ for every $\a \in \mathfrak I,$ then 
$[\mathcal D,\mathcal Y(v,z)] = \frac d{dz} \mathcal Y(v,z)$ for any $v \in \mathcal V$.
\end{prop}

We say that a vertex algebra $\mathcal V$ has a PBW basis, associated with quantum fields $\{X^\a(z)\}_{\a\in\mathfrak I}$, if the index set $\mathfrak I$ is ordered, and we have a linear basis of $\mathcal V$, formed by the vectors
$$\left\{X^{\a_1}_{(-n_1-1)}\dots X^{\a_k}_{(-n_k-1)} \1 \, \biggr| \, 
n_1 \ge n_2 \ge \dots \ge n_k \ge 0, \text { and if } n_i=n_{i+1}, \text { then } \a_i\preceq\a_{i+1}\right\}.$$

For two mutually local quantum fields $X(z),Y(w)$ we introduce the operator product expansion (OPE) formalism, and write
$$X(z)Y(w) \sim \sum_{j} \frac {C_j(w)}{(z-w)^j},$$
if for a finite collection of quantum fields $\{C_j(w)\}_{j=1,2,\dots}$ we have the equality
\begin{equation*}
X(z)Y(w) = \sum_{j} \frac {C_j(w)}{(z-w)^j} \; + :X(z)Y(w):  
\end{equation*}
where $\frac 1{(z-w)^j}$ should be expanded into the Laurent series in non-negative powers of $\frac wz$.
The importance of OPE lies in the fact that all commutators $[X_m,Y_n]$ of Laurent coefficients of quantum fields $X(z),Y(w)$ are
completely encoded by the collection $\{C_j(w)\}$.

\medskip

The remainder of this subsection presents some examples of vertex algebras, which will be used in this paper. All of these algebras are bi-graded and 
have a PBW basis associated with given quantum fields, for which we specify the OPEs.

\begin{example}
We denote by $\hat F(\beta,\gamma)$ the vertex algebra generated by quantum fields
\begin{alignat*}{5}
\beta(z) &= \sum_{n \in \Z} \beta_n  z^{-n}, \qquad &|\beta(z)| &= 0, \qquad \deg \beta(z) &= 0, \\
\gamma(z) &= \sum_{n \in \Z} \gamma_n z^{-n-1}, \qquad &|\gamma(z)| &= 1, \qquad \deg \gamma(z) &= 0,
\end{alignat*}
with the operator product expansions
\begin{equation}\label{eq:beta-gamma OPE} 
\beta(z)\gamma(w) \sim \frac 1{z-w}, \qquad \beta(z) \beta(w) \sim \gamma(z)\gamma(w) \sim 0.
\end{equation}
The commutation relations for the underlying Heisenberg algebra are 
\begin{equation}\label{eq:affine beta-gamma commutation} 
[\beta_m,\gamma_n] = \delta_{m+n,0}, \qquad [\beta_m,\beta_n] = [\gamma_m,\gamma_n]=0.
\end{equation}
\end{example}

\medskip

\begin{example}
We denote by $\hat\Lambda(\psi,\psi^*)$ the vertex superalgebra generated by quantum fields
\begin{alignat*}{5}
\psi(z) &= \sum_{n\in\Z} \psi_n z^{-n-1}, \qquad & |\psi(z)| &= -1, \qquad &\deg \psi(z) &= 1,\\ 
\psi^*(z) &= \sum_{n\in\Z} \psi^*_n z^{-n}, \qquad & |\psi^*(z)| &= 1, \qquad &\deg \psi^*(z) &= 0, 
\end{alignat*}
with the operator product expansions
$$\psi(z)\psi(w) \sim \psi^*(z)\psi^*(w) \sim 0, \qquad \psi(z) \psi^*(w) \sim \frac 1{z-w}.$$
The (anti)-commutation relations for the underlying Clifford algebra are 
\begin{equation}
\label{eq:psi system relations}
\{\psi_m, \psi_n \} = \{ \psi^*_m, \psi^*_n \} = 0, \qquad
\{\psi_m, \psi^*_n\} = \delta_{m+n,0}.
\end{equation}
\end{example}

\medskip

\begin{example}
We denote by $\ghat_k$ the vertex algebra generated by quantum fields
$$X_n = \sum_{n \in \Z} X_n z^{-n-1}, \qquad |X(z)| = 0, \qquad \deg X(z) = 1, \qquad X \in \g,$$
with the operator product expansions
$$X(z)Y(w) \sim \frac {[X,Y](w)}{z-w} + k \, \frac {\<X,Y\>}{(z-w)^2},\qquad k \in \C$$
where $\<\cdot,\cdot\>$ is the Killing form on $\g$. The number $k$ is called the level of $\ghat_k$.
\end{example}
We note that a module for the vertex algebra $\ghat_k$ is a $\Z$-graded $\ghat$-module $\hat V = \bigoplus_{n \ge n_0} \hat V[n]$, such that
$\pi_{\hat V} (\mathbf k) = k \cdot \Id_{\hat V}$ and $\ghat[m] \hat V[n] \subset \hat V[m+n]$ for any $m,n \in \Z$.

\medskip

\begin{example}
We denote by $\vir_c$ the vertex algebra generated by the quantum field
$$L(z) = \sum_{n\in\Z} L_n z^{-n-2}, \qquad  |L(z)| = 0, \qquad \deg L(z) = 2,$$
with the operator product expansion
$$L(z)L(w) \sim \frac {c/2}{(z-w)^4} +  \frac {2 L(w)}{(z-w)^2} + \frac{L'(w)}{z-w}, \qquad c \in \C.$$
The number $c$ is called the central charge of $\vir_c$.
\end{example}
A module for the vertex algebra $\vir_c$ is a $\Z$-graded $\vir$-module $\tilde V = \bigoplus_{n \ge n_0} \tilde V[n]$, such that
$\pi_{\tilde V} (\mathbf c) = c \cdot \Id_{\tilde V}$ and $L_{-m} \tilde V[n] \subset \tilde V[m+n]$ for any $m,n \in \Z$.

\medskip

\begin{example}

We denote by $\hat F_\varkappa(a)$ the vertex algebra generated by the quantum field
$$a(z) = \sum_{n \in \Z} a_n z^{-n-1}, \qquad |a(z)| = 0, \qquad \deg a(z) = 1,$$
with the operator product expansion
\begin{equation}\label{eq:a OPE} 
a(z)a(w) \sim \frac {2\varkappa}{(z-w)^2}, \qquad \varkappa \in \C.
\end{equation}
The commutation relations for the underlying Heisenberg algebra $\mathcal H(a)$ are 
\begin{equation}\label{eq:a commutation}
[a_m, a_n] = 2\varkappa \, m \, \delta_{m+n,0}.
\end{equation}
Note that the operator $a_0$ is central and kills the vacuum.
\end{example}

\medskip

Below we give the construction of a vertex algebra, which will be crucial for our future considerations.
Let $\hat F_{-\varkappa}(\bar a)$ be defined similarly to $\hat F_\varkappa(a)$, so that
\begin{equation}\label{eq:bar a commutation}
[\bar a_m, \bar a_n] = -2\varkappa \, m \,\delta_{m+n,0}, \qquad
\bar a(z) \bar a(w) \sim - \frac {2\varkappa}{(z-w)^2}.
\end{equation}

\begin{thm}
\label{thm:dual Fock vertex structure}
Let $\varkappa \ne 0$. The space $\tilde{\mathbb F}_\varkappa = \hat F_\varkappa(a) \o \hat F_{-\varkappa}(\bar a) \o \C[\mathbf P]$
has a vertex algebra structure, extending those of $\hat F_\varkappa(a)$ and $\hat F_{-\varkappa}(\bar a)$, and such that
$a_0 \1_\lambda = \bar a_0 \1_\lambda = \lambda \1_\lambda$.
\end{thm}

\begin{proof}
Introduce the quantum fields $\{\mathbb Y(\mu,w)\}_{\mu\in\mathbf P}$ by
\begin{equation}\label{eq:double vertex operator}
\begin{split}
\mathbb Y(\mu,z) & = \exp\( \frac {\mu}{2\varkappa}\sum_{n<0} \frac {a_n}{-n} z^{-n}\) 
\exp\(\frac {\mu}{2\varkappa}\sum_{n>0} \frac {a_n}{-n} z^{-n}\) \times\\
& \times
\exp\( -\frac {\mu}{2\varkappa}\sum_{n<0} \frac {\bar a_n}{-n} z^{-n}\) 
\exp\( -\frac {\mu}{2\varkappa}\sum_{n>0} \frac {\bar a_n}{-n} z^{-n}\) \, \1_\mu .
\end{split}
\end{equation}
Straightforward computations lead to the operator product expansions
$$a(z) \bar a(w) \sim \bar a(z) a(w) \sim \mathbb Y(\mu,z) \mathbb Y(\nu,w) \sim 0,$$
$$a(z)\mathbb Y(\mu,w)\sim \bar a(z)\mathbb Y(\mu,w) \sim \frac {\mu \, \mathbb Y(\mu,w) }{z-w},$$
and establish mutual pairwise locality for the quantum fields $a(z),\bar a(z), \mathcal Y(\mu,z)$.

The vacuum is, of course, the vector $\1 \o \1 \o \1_0 \in \mathbb F_\varkappa$. We set
$\mathcal Y(\1_\lambda,z) = \mathbb Y(\lambda,z)$ for any $\lambda \in \mathbf P$.
The spanning and linear independence conditions of 
Proposition \ref{thm:free field construction} are immediate.
Finally, we set $\mathcal D \1_\lambda = \frac \lambda{2\varkappa} (a_{-1} - \bar a_{-1}) \1_\lambda$.
The conditions on $\mathcal D$ amount to
\begin{equation}
\label{eq:double vertex derivative}
\mathbb Y'(\lambda,z) = \frac \lambda{2\varkappa} \,  \biggr(:a(z)\mathbb Y(\lambda,z): - :\bar a(z)\mathbb Y(\lambda,z):\biggr),
\end{equation}
which is checked directly. Applying Proposition \ref{thm:free field construction}, we get the desired statement.
\end{proof}

Theorem \ref{thm:dual Fock vertex structure} should be compared with the construction of lattice vertex algebras.
It is known that the space $\hat F_\varkappa(a) \o \C[\mathbf P]$ carries a vertex algebra structure only for special values of $\varkappa$,
satisfying certain integrality conditions.

\subsection{Bosonic realizations}
We now proceed to study the generalizations of the algebra $\mathfrak R(G_0)$.
As in the finite-dimensional case, we study modules for the Lie algebra $\ghatghat$, which is equivalent to 
having two commuting actions of $\ghat$ on the same space.

As in the classical case, the regular $\ghat$-actions on $\mathfrak R_\varkappa(\hat G_0)$, described in
Theorem \ref{thm:affine action}, can be reformulated in terms of representations of Heisenberg algebras.
We note that the operators
\begin{equation*}
\begin{aligned}
\beta_n &= -y_{-n}\\
\gamma_n &= \pd {y_n}
\end{aligned}, \qquad\qquad
a_n = \begin{cases} 
\pd {\zeta_n}, & n> 0 \\
\zeta \pd {\zeta}, & n = 0 \\
\pd {\zeta_n} - 2 n \varkappa \, \zeta_{-n}, & n < 0
\end{cases}
\end{equation*}
satisfy the commutation relations \eqref{eq:affine beta-gamma commutation},\eqref{eq:a commutation}, and similarly for
\begin{equation*}
\begin{aligned}
\bar\beta_n &= x_{-n},\\
\bar\gamma_n &= - \pd {x_{n}}
\end{aligned}, \qquad\qquad
\bar a_n = \begin{cases} 
\pd {\zeta_n} + 2n \varkappa \, \zeta_{-n}  & n> 0 \\
\zeta \pd {\zeta}, & n = 0\\
\pd {\zeta_{n}} & n < 0\\
\end{cases}.
\end{equation*}
Note also that $\C[\zeta^{\pm1}] \cong \C[\mathbf P]$, and $a_0 = \bar a_0 = a$ act on $\C[\mathbf P]$ by derivations $a \1_\lambda = \lambda \1_\lambda$.

The formulas of Theorem \ref{thm:affine action} are particularly simple, when written for the generating series
$\mathbf e(z), \mathbf h(z), \mathbf f(z)$. For example, \eqref{eq:right affine action} becomes

\begin{equation}
\label{eq:unordered bosons}
\begin{split}
\pi_r(\mathbf e(z)) &= \gamma(z),\\
\pi_r(\mathbf h(z)) &= 2 \beta(z)\gamma(z) + a(z),\\
\pi_r(\mathbf f(z)) &= -\beta(z)^2\gamma(z) - \beta(z)a(z) - \varkappa \beta'(z) +
\exp\(\frac 1\varkappa \sum_{n \ne 0} \frac{a_n- \bar a_n}{n} z^{-n} \) \bar \gamma(z) \1_{-2}.
\end{split}
\end{equation}

Note that in this polynomial realization the constants are annihilated by $\{\mathbf e_n\}_{n \in \Z}$
and $\{\mathbf h_n\}_{n\ge0}$.
The vertex algebra formalism requires a different choice of vacuum, and the introduction
of normal ordering to make products of quantum fields well-defined. This procedure is well-known
in the theory of Wakimoto modules (see \cite{BFr} and references therein), 
for which the $\ghat$-action is constructed by modifying the formulas originating from the
semi-infinite flag variety. In particular, one expects the shifts of the levels of the 
representations by the dual Coxeter number $h^\vee = 2$.

The modifications of the formulas \eqref{eq:unordered bosons} leads to the following result.

\begin{thm}\label{thm:affine bimodule action}
Let $\varkappa \ne 0,$ and let $k = \varkappa - h^\vee$ and $\bar k = -\varkappa - h^\vee$, and let
$$\hat{\mathbb F}_\varkappa = \hat F(\beta,\gamma) \o \hat F(\bar\beta,\bar\gamma) \o \tilde{\mathbb F}_\varkappa.$$
\begin{enumerate}
\item
The space $\hat{\mathbb F}_\varkappa$ has a $\ghatKK$-module structure, defined by
\begin{equation}\label{eq:affine boson left}
\begin{split} 
\mathbf e(z) &= \gamma(z),\\
\mathbf h(z) &= 2:\beta(z) \gamma(z): + a(z),\\
\mathbf f(z) &= -:\beta(z)^2 \gamma(z): - \beta(z) a(z) - k \beta'(z) + \mathbb Y(-2,z) \bar\gamma(z),
\end{split}
\end{equation}
\begin{equation} \label{eq:affine boson right}
\begin{split}
\bar {\mathbf e}(z) &= \bar\gamma(z),\\
\bar {\mathbf h}(z) &= 2:\bar\beta(z) \bar\gamma(z): + \bar a(z),\\
\bar {\mathbf f}(z) &= -:\bar\beta(z)^2 \bar\gamma(z): - \bar\beta(z) \bar a(z) - \bar k \bar\beta'(z) + \mathbb Y(-2,z)\gamma(z).
\end{split}
\end{equation}
\item
The space $\hat {\mathbb F}_\varkappa$ has a compatible VOA structure with $\rank \hat {\mathbb F}_\varkappa = 6$.
(Compatible means that the operators $\mathcal Y(v,z)$ are $\ghatKK$-intertwining operators in the VOA sense).
\end{enumerate}
\end{thm}

Similar formulas for the two commuting actions of $\ghat$ were suggested in \cite{FP}, by analogy with the finite-dimensional
Gauss decomposition of $G$. However, in order to get a meaningful VOA structure - and the corresponding semi-infinite cohomology theory! - 
one must incorporate the twist by $\mathbf {w_0}$, built into the Bruhat decomposition.

One can recover the original Wakimoto realization from \eqref{eq:affine boson left} by properly discarding the 'bar' variables.
We use superscripts 'W' to distinguish the Wakimoto $\ghat_k$-action from \eqref{eq:affine boson left}.

\begin{cor}
The space $\hat W_{\lambda,k} = \hat F(\beta,\gamma)\o \hat F_\varkappa(a) \o \C \1_\lambda$ 
has the structure of a $\ghat_k$-module with $k = \varkappa - h^\vee$, defined by the formulas
\begin{equation}
\begin{split} \label{eq:Wakimoto}
\mathbf e^W(z) &= \gamma(z),\\
\mathbf h^W(z) &= 2:\beta(z) \gamma(z): + a(z),\\
\mathbf f^W(z) &= -:\beta(z)^2 \gamma(z): - \beta(z) a(z) - k\, \beta'(z).
\end{split}
\end{equation}
The $\ghat_k$-module $\hat W_{\lambda,k}$ is called the Wakimoto module.
\end{cor}

\begin{proof}[Proof of Theorem \ref{thm:affine bimodule action}]
It suffices to show that modifying the standard Wakimoto actions by the extra terms 
\begin{align*}
\delta\mathbf f(z) & = \mathbf f(z) - \mathbf f^W(z)  = \mathbb Y(-2,z) \bar\gamma(z), \\
\overline{\delta\mathbf f}(z) & = \bar{\mathbf f}(z) - \bar{\mathbf f}^W(z)  = \mathbb Y(-2,z) \gamma(z),
\end{align*}
does not destroy the operator product expansions.

We begin by showing that the commutation relations for $\ghat_k$ hold. Only those involving the modified quantum field $\mathbf f(z)$ 
need to be considered. We have:
\begin{align*}
\mathbf e^W(z) \, \delta\mathbf f(w) &=  \gamma(z) \mathbb Y(-2,w) \bar\gamma(w) \sim 0,\\
\mathbf h^W(z) \, \delta\mathbf f(w) &= \(2 :\beta(z) \gamma(z): + a(z)\) \mathbb Y(-2,w) \bar\gamma(w) \sim \\
& \sim a(z) \mathbb Y(-2,w) \bar \gamma(w) \sim -\frac {2 \mathbb Y(-2,w)}{z-w} \bar\gamma(w) = - \frac{2 \,\delta\mathbf f(w)}{z-w},\\
\mathbf f^W(z) \, \delta\mathbf f(w) &=  -\beta(z) a(z)\gamma(z) \mathbb Y(-2,w) \bar\gamma(w) \sim \frac {2 \, \mathbb Y(-2,w)}{z-w} \beta(w) \bar\gamma(w),\\
\delta\mathbf f(z) \, \delta\mathbf f(w) &= \mathbb Y(-2,z) \bar\gamma(z) \mathbb Y(-2,w) \bar\gamma(w) \sim 0.
\end{align*}
Using the operator product expansions above we immediately check that
\begin{align*}
\mathbf e(z) \mathbf f(w) &= \mathbf e^W(z) \mathbf f^W(w) + \mathbf e^W(z) \, \delta\mathbf f(w) \sim
\( \frac k{(z-w)^2} + \frac {\mathbf h^W(w)}{z-w} \) + 0 = \frac k{(z-w)^2} + \frac {\mathbf h(w)}{z-w},\\
\mathbf h(z) \mathbf f(w) &= \mathbf h^W(z) \mathbf f^W(w) + \mathbf h^W(z) \, \delta\mathbf f(w) \sim
- \frac {2 \, \mathbf f^W(w)}{z-w} + 0 = - \frac {2 \, \mathbf f(w)}{z-w},\\
\mathbf f(z) \mathbf f(w) &= \mathbf f^W(z) \mathbf f^W(w) + \mathbf f^W(z) \, \delta\mathbf f(w) + 
\delta\mathbf f(z) \, \mathbf f^W(w) + \delta\mathbf f(z) \, \delta\mathbf f(w)  \sim \\
&\sim 0 + \frac {2 \, \mathbb Y(-2,w)}{z-w} \beta(w) \bar\gamma(w) + \frac {2 \, \mathbb Y(-2,z)}{w-z} \beta(z) \bar\gamma(z) + 0 \sim 0,
\end{align*}
and since the operator product expansions not involving $\mathbf f(z)$ are unchanged, we have proved the commutation relations 
for the (left) $\ghat_k$-action. Similarly, one verifies the commutation relations for the (right) $\ghat_{\bar k}$-action.

We now prove that the two actions of $\ghat_k$ and $\ghat_{\bar k}$ commute. We have
\begin{align*}
\bar {\mathbf e}^W(z) \, \delta\mathbf f(w) &=  \bar \gamma(z) \mathbb Y(-2,w) \bar\gamma(w) \sim 0,\\
\bar {\mathbf h}^W(z) \, \delta\mathbf f(w) &= \(2 :\bar \beta(z) \bar \gamma(z): + \bar a(z) \) \mathbb Y(-2,w) \bar\gamma(w) \sim \\
& \sim 2 \frac{ \bar\gamma(z) }{z-w} \mathbb Y(-2,w) - \frac {2 \mathbb Y(-2,w)}{z-w} \bar\gamma(w) \sim 0,
\end{align*}
which implies that 
$\bar{\mathbf e}(z) \mathbf f(w) \sim \bar{\mathbf h}(z) \mathbf f(w) \sim 0.$ Finally, we compute 
\begin{align*}
\bar{\mathbf f}^W(z) \, \delta\mathbf f(w) &= \biggr( -:\bar\beta(z)^2\bar\gamma(z): - \bar k \bar\beta'(z) -
 \bar\beta(z) \bar a(z) \biggr)  \biggr(\mathbb Y(-2,w) \bar\gamma(w) \biggr) \sim \\
& \sim -2 \, \frac{\bar\beta(z) \bar\gamma(z)}{z-w} \mathbb Y(-2,w) + \frac {\bar k}{(z-w)^2} \mathbb Y(-2,w) -\\
& - \( - \frac {2 \mathbb Y(-2,w)}{z-w} :\bar\beta(z) \bar\gamma(w): + 
\frac { :\bar a(w) \mathbb Y(-2,w):}{z-w} - \frac {2 \mathbb Y(-2,w)}{(z-w)^2} \)  \sim \\
& \sim \frac {(\bar k+2) \mathbb Y(-2,w)} {(z-w)^2} - \frac {:\bar a(w)\mathbb Y(-2,w):}{z-w} = 
- \varkappa \frac {\mathbb Y(-2,w)} {(z-w)^2} - \frac {:\bar a(w)\mathbb Y(-2,w):}{z-w}.
\end{align*}
and similarly
\begin{equation*}
\overline{\delta\mathbf f}(z) \mathbf f^W(w) \sim 
\varkappa \frac {\mathbb Y(-2,w)} {(z-w)^2} + \frac {:\bar a(w) \mathbb Y(-2,w):}{z-w} \sim 
- \, \bar{\mathbf f}^W(z) \, \delta\mathbf f(w).
\end{equation*}
Therefore,
\begin{equation*}
\bar{\mathbf f}(z) \mathbf f(w) = 
\bar{\mathbf f}^W(z) \mathbf f^W(w) + \bar{\mathbf f}^W(z) \, \delta\mathbf f(w) + 
\overline{\delta\mathbf f}(z) \, \mathbf f^W(w) +  \overline{\delta\mathbf f}(z) \, \delta\mathbf f(w) \sim 0,
\end{equation*}
and we have established the commutativity of the two $\ghat$-actions.

Proposition \ref{thm:free field construction} implies that $\hat {\mathbb F}_\varkappa$ is a vertex algebra. 
The formulas \eqref{eq:affine boson left} can be written as
\begin{equation*}
\begin{split}
\mathbf e(z) &= \mathcal Y(\gamma_{-1}\1_0,z), \\
{\mathbf h}(z) &= \mathcal Y(2\beta_0\gamma_{-1}\1_0 + a_{-1}\1_0,z), \\
\mathbf f(z) &= \mathcal Y(-(\beta_0)^2\gamma_{-1}\1_0 - a_{-1}\beta_0\1_0 - k \beta_{-1}\1_0 - \bar\gamma_{-1}\1_{-2},z),
\end{split}
\end{equation*}
which means that the quantum fields \eqref{eq:affine boson left} are special cases of the operators $\mathcal Y(\cdot,z)$. 
The same is true for the quantum fields \eqref{eq:affine boson right}.
Therefore, the vertex algebra structure is compatible (in the vertex algebra sense) with the $\ghatKK$-module structure on $\hat{\mathbb F}_\varkappa$.

To give $\hat{\mathbb F}_\varkappa$ a VOA structure we need to introduce the Virasoro element. The Sugawara construction for the affine algebra $\ghat_k$ gives a Virasoro quantum field with central charge 
$c = \frac {3k}{k+h^\vee} = 3 - \frac 6\varkappa$:
\begin{equation}
\label{eq:Sugawara 1}
\begin{split}
L(z) & = \frac 1{2\varkappa} \( \frac{:\mathbf h^2(z):}2 + :\mathbf e(z)\mathbf f(z): + :\mathbf f(z)\mathbf e(z):\) =  \\
& = \frac {:a(z)^2:}{4\varkappa}   - \frac  {a'(z)}{2\varkappa} \, - :\beta'(z) \gamma(z): + 
\frac 1\varkappa \, \mathbb Y(-2,z) \gamma(z) \bar\gamma(z).
\end{split}
\end{equation}
We also note that 
$$L(z) = L^W(z) - \frac 1\varkappa \, \mathbb Y(-2,z) \gamma(z) \bar\gamma(z),$$
where $L^W(z)$ is the Virasoro quantum field given by the Sugawara construction for the standard Wakimoto realization \eqref{eq:Wakimoto}.

Similarly, the affine algebra $\ghat_{\bar k}$ produces another Virasoro quantum field with central charge 
$\bar c = \frac {3\bar k}{\bar k+h^\vee} = 3 + \frac 6\varkappa$:
\begin{equation}
\label{eq:Sugawara 2}
\begin{split}
\bar L(z) & = - \frac 1\varkappa 
\( \frac{:\bar {\mathbf h}^2(z):}2 + :\bar{\mathbf e}(z)\bar{\mathbf f}(z): + :\bar{\mathbf f}(z)\bar{\mathbf e}(z):\) = \\
& = - \frac {:\bar a(z)^2:}{4\varkappa}   + \frac {\bar a'(z)}{2\varkappa} \,  - :\bar \beta'(z) \bar \gamma(z): - 
\frac 1\varkappa \, \mathbb Y(-2,z) \gamma(z) \bar\gamma(z).
\end{split}
\end{equation}

We set $\mathcal L(z) = L(z) + \bar L(z) = L^W(z) + \bar L^W(z)$. To show that $\mathcal L_{-1} = \mathcal D$, we check that 
\begin{equation}\label{eq:affine L_{-1}}
\mathcal Y(\mathcal L_{-1} v,z) = \frac d{dz} \mathcal Y(v,z), \quad v \in \hat{\mathbb F}_\varkappa,
\end{equation}
which for all the generating quantum fields follows from straightforward computations.

Finally, the rank of the VOA $\hat{\mathbb F}_\varkappa$ is equal to 
$$\rank \hat{\mathbb F}_\varkappa = c + \bar c = \(3 - \frac 6\varkappa\) + \(3 + \frac 6\varkappa\) = 6.$$
This concludes the proof of the theorem.
\end{proof}

\subsection{$\ghatKK$-module structure of $\hat{\mathbb F}_\varkappa$ for generic $\varkappa$.}

We now prove the analogue of the Theorem \ref{thm:classical bimodule structure}, 
describing the structure of the $\ghatKK$-module $\hat {\mathbb F}_\varkappa$ for generic values of the parameter $\varkappa$.

For $\lambda \in \h^*, k \in \C$ we denote by $\hat V_{\lambda,k}$ the irreducible $\ghat_k$-module,
generated by a vector $\hat v$ satisfying $\g_+ \hat v = \n_+ \hat v = 0$ and 
$\mathbf h \hat v = \lambda \, \hat v$.

For any $\ghat_k$-module $\hat V$, the restricted dual space $\hat V'$ can be equipped with a $\ghat_k$-action by
$$\<g_n \, v', v\> = - \< v', \hat\omega(g_{-n}) v\>, \qquad v \in \hat V,\  v' \in \hat V', \ g \in \g,$$
where $\hat\omega$ is as in \eqref{eq:affine Cartan automorphism}. We denote the resulting dual module by $\hat V^\star$.

An important source of $\ghat_k$-modules is the induced module construction. Any $\g$-module $V$ may be regarded as a module for the subalgebra $\mathfrak p =  \bigoplus_{n\ge0} \ghat[n]$, with $\g[n]$ acting trivially for $n>0$ and $\mathbf k$ acting as the multiplication by a scalar $k \in \C.$ The induced $\ghat_k$-module $\hat V_k$ is defined as the space
\begin{equation}\label{eq:induced module}
\hat V_k = \mathcal U(\ghat) \o_{\mathcal U(\mathfrak p)} V,
\end{equation}
with $\ghat_k$ acting by left multiplication.

For the remainder of this section, we will assume that complex numbers $\varkappa, k, \bar k$ satisfy
\begin{equation}
\varkappa \notin \mathbb Q, \qquad k = \varkappa - h^\vee, \qquad \bar k = -\varkappa - h^\vee.
\end{equation}

\begin{thm}\label{thm:affine bimodule structure}
There exists a filtration 
\begin{equation}\label{eq:affine filtration}
0 \subset \hat{\mathbb F}_\varkappa^{(0)} \subset \hat{\mathbb F}_\varkappa^{(1)} \subset \hat{\mathbb F}_\varkappa^{(2)} = \hat{\mathbb F}_\varkappa
\end{equation}
of $\ghatKK$-submodules of $\hat{\mathbb F}_\varkappa$ such that
\begin{align}
\hat{\mathbb F}_\varkappa^{(2)}/\hat{\mathbb F}_\varkappa^{(1)} &\cong 
\bigoplus_{\lambda\in\mathbf P^+} \hat V_{-\lambda-2,k} \o \hat V^\star_{-\lambda-2,\bar k}, \label{eq:affine socle F2}\\
\hat{\mathbb F}_\varkappa^{(1)}/\hat{\mathbb F}_\varkappa^{(0)} &\cong 
\bigoplus_{\lambda\in\mathbf P^+} \( \hat V_{\lambda,k} \o \hat V^\star_{-\lambda-2,\bar k} \oplus 
\hat V_{-\lambda-2,k} \o \hat V^\star_{\lambda,\bar k} \), \label{eq:affine socle F1}\\
\hat{\mathbb F}_\varkappa^{(0)} &\cong 
\bigoplus_{\lambda\in\mathbf P} \hat V_{\lambda,k} \o \hat V^\star_{\lambda,\bar k}. \label{eq:affine socle F0}
\end{align}
\end{thm}

\begin{proof}
The operator $\mathcal L_0$ determines a $\Z$-grading $\deg$ of $\hat{\mathbb F}_\varkappa$, which is explicitly described by
\begin{equation}\label{eq:affine grading}
\deg \1_\lambda = 0, \qquad \deg X_n = -n \quad \text{ for } X = a,\bar a, \beta,\gamma, \bar\beta,\bar\gamma.
\end{equation}
The lowest graded subspace $\hat {\mathbb F}_\varkappa[0] = F(\beta_0,\gamma_0) \o F(\bar\beta_0,\bar\gamma_0)\o \C[\mathbf P]$ 
of the vertex algebra $\hat {\mathbb F}_\varkappa$
is identified with the Fock space $\mathbb F$ for the finite-dimensional Lie algebra $\g$. Moreover, since $\varkappa$ is generic, $\hat{\mathbb F}_\varkappa$ 
can be constructed as the induced $\ghatKK$-module from the $\ggbar$-module $\mathbb F$:
$$\hat {\mathbb F}_\varkappa = \mathcal U(\ghatghat)  \o_{\mathcal U(\mathfrak p \oplus \mathfrak p)} \mathbb F.$$

We construct the filtration \eqref{eq:affine filtration} by inducing it from the finite-dimensional one \eqref{eq:classical filtration}:
$$\hat {\mathbb F}_\varkappa^{(0)} =  \mathcal U(\ghatghat)  \o_{\mathcal U(\mathfrak p \oplus \mathfrak p)} \mathbb F^{(0)},
\qquad
\hat {\mathbb F}_\varkappa^{(1)} =  \mathcal U(\ghatghat)  \o_{\mathcal U(\mathfrak p \oplus \mathfrak p)} \mathbb F^{(1)}.
$$

It is easy to check that \eqref{eq:classical socle F2},\eqref{eq:classical socle F1},\eqref{eq:classical socle F0} respectively imply
\eqref{eq:affine socle F2},\eqref{eq:affine socle F1},\eqref{eq:affine socle F0},
which proves the theorem.
\end{proof}

The analogue of the Corollary \ref{thm:classical positive subalgebra}, describing the realization 
of the subalgebra $\mathfrak R(G) \subset \mathfrak R(G_0)$, is given below. 

\begin{thm}\label{thm:affine positive subalgebra}
There exists  a subspace $\hat{\mathbf F}_\varkappa \subset \hat{\mathbb F}_\varkappa$, satisfying
\begin{enumerate}
\item $\hat{\mathbf F}_\varkappa$ is a vertex operator subalgebra of $\hat{\mathbb F}_\varkappa$, and is generated 
by the quantum fields \eqref{eq:affine boson left}, \eqref{eq:affine boson right} and $\mathbb Y(1,z)$. 
In particular, $\mathbf F$ is a $\ghatKK$-submodule of $\hat{\mathbb F}_\varkappa$.
\item As a $\ghatKK$-module, $\hat{\mathbf F}_\varkappa$ is generated by the vectors $\{\1_\lambda\}_{\lambda \in \mathbf P^+}$, and we have
\begin{equation}
\label{eq:affine positive decomposition}
\hat{\mathbf F}_\varkappa \cong \bigoplus_{\lambda \in \mathbf P^+} \hat V_{\lambda,k} \o \hat V^\star_{\lambda, \bar k}.
\end{equation}
\end{enumerate}
\end{thm}

\begin{proof}
As before, we identify the lowest graded subspace $\hat {\mathbb F}_\varkappa[0]\subset \hat {\mathbb F}_\varkappa$ with the Fock space 
$\mathbb F$ for the finite-dimensional Lie algebra $\g$. Recall from Corollary \ref{thm:classical positive subalgebra} that the 
$\ggbar$-module $\mathbb F$ contains the distinguished submodule $\mathbf F$. We define the subspace $\hat {\mathbf F}_\varkappa$ 
as the $\ghatKK$-submodule of $\hat {\mathbb F}_\varkappa$, induced from $\mathbf F$:
$$\hat {\mathbf F}_\varkappa = \mathcal U(\ghatghat)  \o_{\mathcal U(\mathfrak p \oplus \mathfrak p)} \mathbf F.$$
It immediately follows from Corollary \ref{thm:classical positive subalgebra} that $\hat{\mathbf F}_\varkappa$ 
is generated by the vectors $\{\1_\lambda\}_{\lambda \in \mathbf P^+}$, and has the decomposition \eqref{eq:affine positive decomposition}.

Next, we need to show that $\hat {\mathbf F}_\varkappa$ is a vertex subalgebra. Let $\hat{\mathbf F}_\varkappa'$ denote 
the space, spanned by the Laurent coefficients of $\hat{\mathbb F}_\varkappa$-valued fields $\mathcal Y(a,z)b$ for all possible 
$a,b \in \hat {\mathbf F}_\varkappa$. We will establish that $\hat{\mathbf F}_\varkappa' = \hat{\mathbf F}_\varkappa$.

Indeed, $\hat{\mathbf F}_\varkappa'$ is a $\ghatKK$-submodule of $\hat{\mathbb F}_\varkappa$, and can be induced from its lowest graded component
$\mathbf F' = \hat{\mathbf F}_\varkappa'[0]$, which is a $\ggbar$-submodule of $\mathbb F$. It suffices to prove that $\mathbf F' = \mathbf F$.

For any $a,b \in \mathbf F$, the lowest graded component of $\mathcal Y(a,z)b$ is equal to the product $a b$ in the algebra $\mathbb F$, 
and since $\mathbf F$ is a subalgebra, we have $a b \in \mathbf F$. (Note that any element $a \in \mathbf F$ can be obtained this way, 
for example, by taking $b = \1$). Using the commutation relations with the two copies of $\ghat$, 
we can prove that the lowest graded component of $\mathcal Y(a,z)b$ lies in $\mathbf F$ for any $a,b \in \hat {\mathbf F}_\varkappa$.

It follows that $\mathbf F' = \mathbf F$ and hence $\hat{\mathbf F}_\varkappa' = \hat{\mathbf F}_\varkappa$, which means that the restrictions 
of the operators $\mathcal Y(\cdot,z)$, corresponding to the subspace $\hat{\mathbf F}_\varkappa$, are well-defined. 
Thus $\hat{\mathbf F}_\varkappa$ is a vertex subalgebra of $\hat{\mathbb F}_\varkappa$. It is clear that as a vertex subalgebra 
$\hat{\mathbf F}_\varkappa$ is generated by the quantum fields \eqref{eq:affine boson left}, \eqref{eq:affine boson right} and 
$\{\mathbb Y(\lambda,z)\}_{\lambda \in \mathbf P^+}$, and the latter are generated by the single operator $\mathbb Y(1,z)$.

Finally, $\hat {\mathbf F}_\varkappa$ contains both $L^W(z)$ and $\bar L^W(z)$ - hence also $\mathcal L(z)$ - and
therefore is a vertex operator subalgebra of $\hat {\mathbb F}_\varkappa$.
\end{proof}

\begin{rem}
The vertex operator algebras $\hat{\mathbf F}_\varkappa$ and $\hat{\mathbb F}_\varkappa$ give explicit realizations of the
modified regular representations $\mathfrak R'_\varkappa(\hat G)$ and $\mathfrak R'_\varkappa(\hat G_0)$ we discussed in the introduction. 
It would be interesting to construct them invariantly by using the correlation functions approach \cite{FZ}, interpreting the
rational functions $\<\1', \mathcal Y(v_1,z_1) \dots \mathcal Y(v_n,z_n) \1\>$ for $v_1,\dots,v_n \in \mathbb F \cong \hat{\mathbb F}_\varkappa[0]$ 
as solutions of differential equations similar to the Knizhnik-Zamolodchikov equations.
\end{rem}

\subsection{Semi-infinite cohomology of $\ghat$}

The fact that the level of the diagonal action of $\ghat$ in the modified regular representations is equal to the special value 
$-2 h^\vee$ allows us to introduce the semi-infinite cohomology of $\ghat$ with coefficients in $\hat {\mathbb F}_\varkappa$ 
and in $\hat {\mathbf F}_\varkappa$. In this section we show that for generic values of $\varkappa$ these cohomologies 
lead to the same algebras of formal characters as in the finite-dimensional case.

We recall the definition of the semi-infinite cohomology \cite{Fe,FGZ}. The main new ingredient is the "space of semi-infinite forms" 
$\boldsymbol{\hat\Lambda}^\semiinfty$, which replaces the finite-dimensional exterior algebra $\boldsymbol\Lambda$. 
We summarize its properties in the following

\begin{prop}
Let $\boldsymbol{\hat\Lambda}^\semiinfty = \bigwedge \ghat_- \o \bigwedge (\ghat'_+ \oplus \g')$. Then

\begin{enumerate}
\item
The Clifford algebra, generated by $\{\iota(g_n), \eps(g'_n)\}_{g \in \g, g' \in \g',n \in \Z}$ with relations
\begin{equation}
\label{eq:affine Clifford relations}
\{\iota(x_m), \iota(y_n) \} = \{ \eps(x'_m), \eps(y'_n) \} = 0, \qquad
\{\iota(x_m), \eps(y'_n)\} = \delta_{m,n}\, \<y',x\>.
\end{equation}
acts irreducibly on $\boldsymbol{\hat\Lambda}^\semiinfty$, so that for any 
$\omega_- \in \bigwedge \ghat_-,\  \omega_+ \in \bigwedge (\ghat'_+ \oplus \g')$ we have
\begin{align*}
\iota(x_n) ( \omega_- \o 1 ) &= 
 \begin{cases} 0,& n \ge 0 \\ (x_n \wedge \omega_-) \o 1, & n<0 \end{cases}, \qquad
\eps(x'_n) (1 \o \omega_+) &=  
\begin{cases}  1 \o (x'_n \wedge \omega_+), & n\ge0\\ 0,& n<0  \end{cases}.
\end{align*}

\item $\boldsymbol{\hat\Lambda}^{\semiinfty}$ is a bi-graded vertex superalgebra, with vacuum $\1 = 1 \o 1$,
and generated by
\begin{alignat*}{5}
\iota(x,z) &= \sum_{n\in\Z} \iota(x_n) z^{-n-1},  \qquad  | \iota(x,z) | & = & -1, \quad &\deg \iota(x,z) &= 1, \quad & x &\in \g,\\
\eps(x',z) &= \sum_{n\in\Z} \eps(x'_{-n}) z^{-n}, \qquad  | \eps(x',z) | & = & 1, \quad &\deg \eps(x',z) &= 0, \quad  &x' &\in \g'.
\end{alignat*}

\item $\boldsymbol{\hat\Lambda}^\semiinfty$ has a $\ghat$-module structure on the level $\mathbf k = 2 h^\vee$, defined by
$$\pi(x_n) = \sum_{m\in\Z} \sum_i :\eps((g'_i)_m)\iota([g_i,x]_{n+m}):,
 \qquad x \in \g.$$

\end{enumerate}
\end{prop}

One can think of $\bigwedge \ghat_-$ as the space spanned by formal ``semi-infinite'' forms
$$\omega = \xi'_{i_1} \wedge \xi'_{i_2} \wedge \xi'_{i_3} \wedge \dots,\qquad i_{n+1} = i_n + 1 \text{ for } n\gg0,$$
where $\{\xi_j\}_{j\in \mathbb N}$ is a homogeneous basis of $\ghat_-$. A monomial 
$\xi_{j_1} \wedge \dots \wedge \xi_{j_m} \in \bigwedge \ghat_-$ is identified with the semi-infinite form with the corresponding factors missing: 
$$\omega = \pm \ \xi'_1 \wedge \xi'_2 \wedge \dots \wedge \xi'_{j_1-1} \wedge \widehat{\xi'_{j_1}} \wedge  \xi'_{j_1+1} \wedge \dots \wedge \xi_{j_m-1} \wedge \widehat{\xi_{j_m}} \wedge \xi_{j_m+1} \wedge\dots.$$
In other words, if $\xi_j \in \ghat$, then $\iota(\xi_j)$ operates as usual by eliminating the factor $\xi'_j$.

\begin{defn}
The BRST complex, associated with a $\ghat$-module $\hat V$ on the level $k = -2h^\vee$, is the complex
$C^{\semiinfty+\bullet} (\ghat, \C\mathbf k;\hat V) = \boldsymbol{\hat\Lambda}^{\semiinfty+\bullet} \otimes \hat V$, with the differential
\begin{equation}\label{eq:affine differential}
\hat{\mathbf d} = \sum_{n\in\Z} \sum_i \eps((g'_i)_n) \pi_{\hat V}((g_i)_n) - \frac12 \sum_{m,n\in \Z} \sum_{i,j} :\eps((g_i)'_m) \eps((g_j)'_n) \iota([g_i,g_j]_{m+n}):,
\end{equation}
where $\{g_i\}$ is any basis of $\g$, and $\{g'_i\}$ is the dual basis of $\g'$.
The corresponding cohomology is denoted $H^{\semiinfty+\bullet}(\ghat,\C \mathbf k;\hat V)$.
\end{defn}

The BRST complex above gives the relative (to the center) version of the semi-infinite cohomology.
We don't consider any other type of cohomology, and thus simply drop the word 'relative' everywhere. The condition $k = -2h^\vee$ 
is equivalent to $\hat{\mathbf d}^2 = 0$.

If $\hat V$ is a vertex algebra, then its semi-infinite cohomology inherits a vertex superalgebra structure \cite{LZ}.

The following theorem is similar to the reduction theorem of \cite{FGZ} (see also \cite{Li}), and relates the semi-infinite cohomology for generic values of $\varkappa$ with the classical cohomology of Lie algebras.

\begin{thm}\label{thm:reduction theorem}
Let $V$ be a $\ggbar$-module, and let $\varkappa \in \C$ be generic. Set $k = \varkappa - h^\vee$ and $\bar k = -\varkappa - h^\vee$, and denote $\hat V$ be the induced $\ghat_k \oplus \ghat_{\bar k}$-module.
Then with respect to the diagonal $\g$-action $\hat V$ is a level $\mathbf k = -2 h^\vee$ module, and
$$H^{\semiinfty+\bullet}(\ghat,\C\mathbf k; \hat V) \cong H^\bullet(\g,V).$$
\end{thm}

\begin{proof}
As a vector space, the module $\hat V$ has a decomposition
$$\hat V = \mathcal U(\ghat_-) \o V \o \mathcal U(\ghat_+)',$$
where we identified the factor $\mathcal U(\ghat_-)$, coming from the right induced action of $\ghat_{\bar k}$, with $\mathcal U(\ghat_+)'$ using the 
non-degenerate (since $\varkappa$ is generic!) contravariant pairing.
Therefore, as vector spaces 
\begin{equation}
\label{eq:affine complex factorization}
C^\bullet(\ghat, \C\mathbf k; \hat V) = C_-^\bullet \o C_0^\bullet \o C_+^\bullet,
\end{equation}
where
$$C_-^\bullet = \bigwedge \ghat_- \o \mathcal U(\ghat_-), \qquad
C_0^\bullet = \bigwedge \g' \o V, \qquad
C_+^\bullet = \bigwedge \ghat'_+ \o \mathcal U(\ghat_+)^*$$

We write the differential $\hat{\mathbf d}$ as
$$\hat{\mathbf d} = \mathbf d_- +  \mathbf d_0 + \mathbf d_+ + \boldsymbol\delta,$$
where $\mathbf d_\pm$ are the BRST differentials for $\ghat_\pm$,
\begin{align*}
\mathbf d_- &=  \sum_{n<0} \sum_i \eps((g'_i)_n) \pi_l((g_i)_n) - \frac 12 \sum_{m,n<0} \sum_{i,j} :\eps((g_i)'_m) \eps((g_j)'_n) \iota([g_i,g_j]_{m+n}):,\\
\mathbf d_+ &= \sum_{n>0} \sum_i \eps((g'_i)_n) \pi_r((g_i)_n) - \frac12 \sum_{m,n>0} \sum_{i,j} :\eps((g_i)'_m) \eps((g_j)'_n) \iota([g_i,g_j]_{m+n}): ,
\end{align*}
the differential $\mathbf d_0$ is defined as in \eqref{eq:finite differential} with the $\g$-action $\pi_V$ replaced by
\begin{equation}
\label{eq:spectral g-action} 
\pi(x) = \pi_{\hat V}(x) + \sum_{n\ne0}  :\eps((g_j)'_n) \iota([x,g_j]_n):, \qquad x \in \g,
\end{equation}
and $\boldsymbol\delta$ includes all the remaining terms:
\begin{align*}
\boldsymbol \delta &= \sum_{n>0} \sum_i \eps((g'_i)_n) \pi_l((g_i)_n) 
+ \sum_{n<0} \sum_i \eps((g'_i)_n) \pi_r((g_i)_n) - \\
&- \sum_{m>0,n<0} \sum_{i,j} :\eps((g_i)'_m) \eps((g_j)'_n) \iota([g_i,g_j]_{m+n}):.
\end{align*}

Following \cite{FGZ}, we introduce the skewed degree $f\deg$ by
$$f\deg (w_- \o w_0 \o w_+) = \deg w_+ - \deg w_-, \qquad w_\pm \in C_\pm,\ \  w_0 \in C_0,$$
where the '$\deg$' gradings in the complexes $C_\pm$ are inherited from 
$C^{\semiinfty}(\ghat,\mathbf k;\hat V)$. We set
$$\mathfrak B^p = \left\{ v \in C^{\semiinfty}(\ghat,\mathbf k;\hat V) \ \biggr| \ f\deg v \ge p \right\}.$$
One can check that $\mathbf d_\pm$ and $\mathbf d_0$ preserve the filtered degree, 
and that $\boldsymbol \delta (\mathfrak B^p) \subset \mathfrak B^{p+1}$. 
Thus, $\{\mathfrak B^p\}_{p\in\Z}$ is a decreasing filtration of the complex 
$C^{\semiinfty}(\ghat,\mathbf k;\hat V)$,
and the associated graded complex has the reduced differential
$$\mathbf d_{red} = \mathbf d_- + \mathbf d_0 + \mathbf d_+.$$
We now compute the corresponding reduced cohomology, which will provide a bridge to 
$H^{\semiinfty+\bullet}(\ghat,\C\mathbf k;\hat V)$.

It is clear that $\mathbf d_\pm^2 = \{\mathbf d_+,\mathbf d_-\} = 0$, and that the differentials $\mathbf d_\pm: C_\pm^\bullet \to C_\pm^{\bullet+1}$ act in their respective factors of \eqref{eq:affine complex factorization}. One can also check that $(\mathbf d_0)^2 = \{\mathbf d_0,\mathbf d_\pm\} = 0$; moreover, 
$$\mathbf d_0 (C_-^\bullet \o C_0^\bullet \o C_+^\bullet) \subset  (C_-^\bullet \o C_0^{\bullet+1} \o C_+^\bullet)$$
despite the fact that $\mathbf d_0$ does not act in $C_0^\bullet$.

It is a well-known fact in homological algebra that
$$H^n(C_+,\mathbf d_+) = H^n(\ghat_+;\mathcal U(\ghat_+)') = \delta_{n,0} \, \C,$$
with $1 \o 1' \in C_+$ representing the non-trivial cohomology class. Similarly, one has 
$$H^n(C_-,\mathbf d_-) = H_{-n}(\ghat_-;\mathcal U(\ghat_-)) = \delta_{n,0} \, \C,$$
and $1 \o 1 \in C_-$ represents the non-trivial cohomology. Further, one can check that
the subspace $1 \o C_0^\bullet \o 1 \subset C^{\semiinfty+\bullet}(\ghat,\C\mathbf k;\hat V)$ is stabilized by $\mathbf d_0$, and that the $\g$-action \eqref{eq:spectral g-action} on that subspace reduces to the $\g$ action $1 \o \pi_V \o 1$.
It follows that
$$H^\bullet_{red}(\ghat, \C\mathbf k; \hat V) \cong H^\bullet(1 \o  C_0 \o 1,\mathbf d_0) \cong H^\bullet(C_0,\mathbf d) \cong H^\bullet(\g,V).$$

We now return to the cohomology of $C^{\semiinfty+\bullet}(\ghat,\C\mathbf k; \hat V)$. Since $\hat {\mathbf d}$ preserves the '$\deg$' grading, it can be computed separately for each subcomplex $C^{\semiinfty+\bullet}(\ghat,\C\mathbf k; \hat V)[m], m \in \Z$.
The filtration $\{\mathfrak B^p[m]\}_{p \in \Z}$ of this complex is finite for each $m$, and leads to a  finitely converging spectral sequence with $E_1^{p,q}[m] = H^q_{red}(\mathfrak B^p[m]/ \mathfrak B^{p+1}[m])$.

For $m\ne 0$ we have $H^q_{red}(\mathfrak B^p[m]/ \mathfrak B^{p+1}[m])=0$ for all $p$, hence the spectral sequence is zero, and  $H^{\semiinfty+\bullet}(\ghat,\C\mathbf k; \hat V)[m]=0$. For $m=0$ we note that
$$\mathfrak B^0[0] = C^{\semiinfty+\bullet}(\ghat,\C\mathbf k; \hat V)[0], \qquad 
\mathfrak B^1[0] = 0,$$
which means that $E_1^{p,q}[0] = 0$ unless $p=0$, and the collapsing spectral sequence implies
$$H^{\semiinfty+\bullet}(\ghat,\C\mathbf k; \hat V)[0] \cong H^\bullet_{red}(\mathfrak B^0[0]/ \mathfrak B^1[0]) \cong H^\bullet(\g;V).$$

This completes the proof of the theorem.
\end{proof}

\begin{cor}
\label{thm:affine cohomology}
The vertex superalgebras $H^{\semiinfty+\bullet}(\ghat,\C\mathbf k; \hat{\mathbf F}_\varkappa), 
H^{\semiinfty+\bullet}(\ghat,\C\mathbf k; \hat{\mathbb F}_\varkappa)$ degenerate into commutative
superalgebras. Moreover, we have commutative superalgebra isomorphisms
\begin{equation}
\label{eq:affine semiinfinite isomorphisms}
H^{\semiinfty+\bullet}(\ghat,\C\mathbf k; \hat{\mathbf F}_\varkappa) \cong H^\bullet(\g; \mathbf F),\qquad
H^{\semiinfty+\bullet}(\ghat,\C\mathbf k; \hat{\mathbb F}_\varkappa) \cong H^\bullet(\g; \mathbb F).
\end{equation}
In particular,
$$H^{\semiinfty+0}(\ghat,\C \mathbf k;\hat{\mathbf F}_\varkappa) \cong 
H^{\semiinfty+0}(\ghat,\C \mathbf k;\hat{\mathbb F}_\varkappa) \cong \C[\mathbf P]^W.$$
\end{cor}

\begin{proof}
Theorem \ref{thm:reduction theorem} gives us isomorphisms \eqref{eq:affine semiinfinite isomorphisms} on the level of vector spaces. 
It is also clear from its proof that the semi-infinite cohomology is concentrated in the subspace of $\deg = 0$, 
and thus the operators $\mathcal Y(\cdot,z)$ on cohomology are reduced to their constant terms. 
In particular, they are independent of $z$, which means that the vertex superalgebra degenerates into a commutative superalgebra. 
The multiplication is easily traced back to the multiplications in $\mathbf F \cong \hat{\mathbf F}_\varkappa[0]$ and in the exterior algebra 
$\boldsymbol\Lambda = \bigwedge \g'$, which shows that \eqref{eq:affine semiinfinite isomorphisms} are superalgebra isomorphisms.
\end{proof}

\section{Modified regular representations of the Virasoro algebra.}

\subsection{Virasoro algebra and the quantum Drinfeld-Sokolov reduction}
In this section we present a construction of the regular representation of the Virasoro algebra, which goes in parallel with
constructions in the previous sections. However, instead of beginning with a space of functions on the corresponding group
(which is, strictly speaking, a semigroup in the complex case), we will use the quantum Drinfeld-Sokolov reduction \cite{FeFrDS} (see also \cite{BFr} and references therein), applied to the modified regular representations of $\ghat$ constructed in Section 1. As a result we obtain certain bimodules over the Virasoro algebra, which have the structure similar to their affine counterparts. The result of the quantum Drinfeld-Sokolov reduction applied to the actual regular representation of $\ghat$ should have a standard interpretation in terms of the space of functions on the Virasoro semigroup, but we will not need this fact for our purposes.

Recall that the Virasoro algebra $\vir$ is the infinite-dimensional complex Lie algebra, generated by $\{L_n\}_{n \in \Z}$ and a central element $\mathbf c,$ subject to the commutation relations
$$[L_m, L_n] = (m-n) L_{m+n} + \frac {m^3-m}{12}\, \mathbf c.$$

The Virasoro algebra has a $\Z$-grading $\vir = \oplus_{n\in\Z} \vir[n]$, determined by 
$$\deg L_n = -n, \qquad \deg \mathbf c = 0.$$

There is a functorial correspondence between certain representations of affine Lie algebras and their $\mathcal W$-algebra counterparts,
called the quantum Drinfeld-Sokolov reduction \cite{FeFrDS} (see also \cite{BFr} and references therein). 
We review this procedure for the case $\ghat = \slhat$,
when the corresponding $\mathcal W$-algebra is identified with the Virasoro algebra.

\begin{defn} For any $\ghat_k$-module $\hat V$, the complex $(C_{DS}(\hat V),\mathbf d_{DS})$,
$$C_{DS}(\hat V) = \hat V \o \hat\Lambda(\psi,\psi^*), \qquad \mathbf d_{DS} = \sum_{n\in\Z} \psi^*_n \pi_{\hat V}(\mathbf e_n) + \psi^*_1,$$
is called the BRST complex of the quantum Drindeld-Sokolov reduction. The corresponding cohomology is denoted $H_{DS}(\hat V)$.
\end{defn}

The BRST complex above is very similar to the semi-infinite cohomology complex for the nilpotent loop algebra
$\nhat_+ = \bigoplus_{n\in\Z} \C \mathbf e_n$.
Indeed, the corresponding space of semi-infinite forms $\boldsymbol{\Lambda}^\semiinfty(\nhat_+)$ is identified 
with $\hat\Lambda(\psi,\psi^*)$ by $\iota(\mathbf e_n) \equiv \psi_n, \ \eps(\mathbf e'_n) \equiv \psi^*_{-n}$,
and the only modification is the additional term $\psi^*_1$ in the differential.

The BRST complex inherits the gradings $|\cdot|$ and $\deg$ from $\hat\Lambda(\psi,\psi^*)$ and $\hat V$.
Since $|\mathbf d_{DS}| = 1$, the grading $|\cdot|$ descends to the cohomology $H_{DS}(\hat V)$.
However, with respect to the other grading, the differential $\mathbf d_{DS}$ is not homogeneous.
We introduce a modified grading $\deg'$ by
\begin{gather*}
\deg' \mathbf e(z) = 0, \qquad
\deg' \mathbf h(z) = 1,\qquad
\deg' \mathbf f(z) = 2, \\
\deg' \psi(z) = 0,\qquad
\deg' \psi^*(z) = 1.
\end{gather*}
The differential $\mathbf d_{DS}$ then satisfies $\deg' \mathbf d_{DS} = 0$, and the grading $\deg'$ descends to $H_{DS}(\hat V)$.

\medskip

The cohomology $H^0_{DS}(\ghat_k)$ of the vacuum module inherits a vertex algebra structure. We have the following result
(details of the proof can be found in \cite{BFr}).

\begin{prop} \label{thm:Drinfeld-Sokolov}
For $k \ne -h^\vee$ we have $H^0_{DS}(\ghat_k) \cong \vir_c$, where
$c = 1 - \frac{6}{k+h^\vee} - 6k$.
\end{prop}

For any $\ghat_k$-module $\hat V$, the vertex algebra $\vir_c \cong H^0_{DS}(\ghat_k)$ acts on $H^0_{DS}(\hat V)$.
For $\varkappa \ne 0$, set
$\tilde F_{\lambda,\varkappa} = H_{DS}^0(\hat W_{\lambda,\varkappa-h^\vee}).$
The following identifies the $\vir_c$-module structure on $\tilde F_{\lambda,\varkappa}$.

\begin{prop}
\label{thm:Wakimoto to FF}
Let $\varkappa \ne 0$, and let $c = 13 - 6 \varkappa - \frac6\varkappa$. Then 
$\tilde F_{\lambda,\varkappa} \cong \hat F_\varkappa(a) \o \C \1_\lambda$
as a vector space, and the $\vir_c$-action is given by
\begin{equation}\label{eq:Feigin-Fuks}
L^F(z) = \frac 1{4\, \varkappa} :a(z)^2: + \frac{\varkappa-1}{2\, \varkappa} \, a'(z).
\end{equation}
\end{prop}

\begin{proof}
In the vector space factorization of $\hat W_{\lambda,\varkappa-h^\vee} = \hat F(\beta,\gamma) \o \hat F_\varkappa(a) \otimes \C \1_\lambda$,
the differential $\mathbf d_{DS}$ acts only in the first component. Therefore, we must have
$$\tilde F_{\lambda,\varkappa} = H_{DS}(\hat F(\beta,\gamma)) \o \hat F_\varkappa(a) \otimes \C \1_\lambda.$$
A spectral sequence reduces the cohomology $H_{DS}^0(\hat F(\beta,\gamma))$ to the cohomology of
the semi-infinite Weil complex $\hat F(\beta,\gamma)\o\hat\Lambda(\psi,\psi^*)$. The latter splits into an infinite product of finite-dimensional Weil
complexes, and thus has one-dimensional cohomology, concentrated in degree $0$.

The inclusion of vertex algebras $\ghat_{\varkappa-h^\vee} \hookrightarrow \hat W_{0,\varkappa-h^\vee}$ induces an inclusion 
$\vir_c \hookrightarrow \tilde F_{0,\varkappa}$, and the explicit formula \eqref{eq:Feigin-Fuks} for $L(z)$ in terms of 
$a(z)$ is a result of a direct computation.
\end{proof}

The realization \eqref{eq:Feigin-Fuks} of Virasoro modules was known long before the quantum Drinfeld-Sokolov reduction, and
is called the Feigin-Fuks construction in the literature.
We use the superscript "F" to distinguish this standard action from the modified Virasoro actions, which we will be considering later.

\subsection{Bosonic realization of the regular representation}

The Virasoro analogue of the Peter-Weyl theorem is more subtle than in the case of classical and affine Lie algebras. 
There is no clear way to calculate the two commuting $\vir$-actions in a way similar to 
Theorem \ref{thm:classical bimodule action} and Theorem \ref{thm:affine action}. 
However, there exists a Fock space realization analogous to Theorem \ref{thm:affine bimodule action}, which we will call
the regular representation of the Virasoro algebra.

\begin{thm}\label{thm:Virasoro action}
Let $\varkappa \ne 0$, and let $c = 13 - 6 \varkappa - \frac 6\varkappa$ and 
$\bar c = 13 + 6 \varkappa + \frac 6\varkappa$.
\begin{enumerate}
\item
The space $\tilde {\mathbb F}_\varkappa$ has a $\virCC$-module structure, defined by
\begin{align}
L(z) &= \frac 1{4\varkappa} :a(z)^2: + \frac{\varkappa-1}{2\varkappa} a'(z)
- \frac 1\varkappa \, \mathbb Y(-2,z), \label{eq:Virasoro action 1}\\
\bar L(z) &= -\frac 1{4\varkappa} :\bar a(z)^2: + \frac{\varkappa+1}{2\varkappa} \bar a'(z)
 + \frac 1\varkappa \, \mathbb Y(-2,z). \label{eq:Virasoro action 2}
\end{align}
\item
The space $\tilde {\mathbb F}_\varkappa$ has a compatible VOA structure with $\rank \tilde {\mathbb F}_\varkappa = 26$.
\end{enumerate}
\end{thm}

\begin{proof}
The formulas \eqref{eq:Virasoro action 1},\eqref{eq:Virasoro action 2} are nothing else but the result of the two-sided quantum Drinfeld-Sokolov reduction,
which consists of two reductions applied separately to the two commuting $\ghat$-actions of Theorem \ref{thm:affine bimodule action}, cf. formulas \eqref{eq:Sugawara 1},\eqref{eq:Sugawara 2} and Proposition \ref{thm:Wakimoto to FF}.

Rather than give detailed proof of this fact, we choose to verify the commutation relations directly. Introduce notation 
\begin{align*}
\delta L(z) &= L(z) - L^F(z) = \frac 1\varkappa \, \mathbb Y(-2,z), \\
\overline{\delta L}(z) &= \bar L(z) - \bar L^F(z) = - \frac 1\varkappa \, \mathbb Y(-2,z).
\end{align*}
Without the additional terms $\delta L(z), \overline{\delta L}(z)$, both \eqref{eq:Virasoro action 1} and
\eqref{eq:Virasoro action 2} give two commuting copies of the standard construction \eqref{eq:Feigin-Fuks} with the specified central charges.
Therefore, it suffices to show that the presence of these extra terms does not violate the commutation relations for
$\virCC$.

Straightforward computations immediately show that
$$\delta L(z) \, \delta L(w) \sim \delta L(z) \, \overline{\delta L}(w)  \sim \overline{\delta L}(z) \, \overline{\delta L}(w) \sim 0.$$

$$L^F(z) \mathbb Y(-2,w) \sim \frac {\mathbb Y(-2,w)}{(z-w)^2} - \frac 1{\varkappa} \frac{:a(w) \mathbb Y(-2,w):}{z-w},$$
$$\bar L^F(z) \mathbb Y(-2,w) \sim \frac {\mathbb Y(-2,w)}{(z-w)^2} + \frac 1{\varkappa} \frac{:\bar a(w) \mathbb Y(-2,w):}{z-w}.$$

We now prove the commutation relations for the action \eqref{eq:Virasoro action 1}. We have
\begin{align*}
L(z)L(w) & - L^F(z) L^F(w)  =  L^F(z) \, \delta L(w) + \delta L(z) L^F(w) + \delta L(z) \, \delta L(w) \sim\\
& \sim   \frac 1\varkappa \( \frac {\mathbb Y(-2,w)}{(z-w)^2} 
- \frac 1{\varkappa} \frac{:a(w) \mathbb Y(-2,w):}{z-w}\)
+ \frac 1\varkappa \(\frac {\mathbb Y(-2,z)}{(z-w)^2} 
+ \frac 1{\varkappa} \frac{:a(z) \mathbb Y(-2,w):}{z-w} \) \sim\\
& \sim  \frac 1\varkappa \( \frac {2 \, \mathbb Y(-2,w)}{(z-w)^2} + \frac {\mathbb Y'(-2,w)}{z-w} \) \sim
\frac {2 \, \delta L(w)}{(z-w)^2} + \frac {(\delta L)'(w)}{z-w} ,
\end{align*}
and therefore
\begin{align*}
L(z)L(w) & \sim L^F(z) L^F(w)  + \frac {2 \, \delta L(w)}{(z-w)^2} + \frac {(\delta L)'(w)}{z-w} \sim
 \( \frac{c/2}{(z-w)^4} + \frac {2L^F(w)}{(z-w)^2} + \frac {(L^F)'(w)}{z-w} \) + \\
& + \frac {2 \, \delta L(w)}{(z-w)^2} + \frac {(\delta L)'(w)}{(z-w)^2} =
\frac{c/2}{(z-w)^4} + \frac {2L(w)}{(z-w)^2} + \frac {L'(w)}{z-w} .
\end{align*}
We have established that adding the extra term $\delta L(z)$ to the action \eqref{eq:Feigin-Fuks} preserves the commutation relations 
for $\vir_c.$ Similarly, the formula \eqref{eq:Virasoro action 2} gives a representation of $\vir_{\bar c}.$

We now show that the two actions of $\vir_c$ and $\vir_{\bar c}$ commute. Using \eqref{eq:double vertex derivative}, we get 
\begin{equation*}
\begin{split}
\delta L(z) \bar L^F(w) & \sim  \frac 1\varkappa \( \frac {\mathbb Y(-2,z)}{(z-w)^2} - \frac 1{\varkappa} \frac{:\bar a(z) \mathbb Y(-2,z):}{z-w} \) \sim \\
&  \sim  \frac 1\varkappa \( \frac {\mathbb Y(-2,w)}{(z-w)^2} + \frac {\mathbb Y'(-2,w)}{z-w} - 
\frac 1{\varkappa} \frac{:\bar a(w) \mathbb Y(-2,w):}{z-w} \) \sim \\
&  \sim  \frac 1\varkappa \( \frac {\mathbb Y(-2,w)}{(z-w)^2} - \frac 1{\varkappa} \frac{:a(w) \mathbb Y(-2,w):}{z-w} \) .
\end{split}
\end{equation*}
Note that \eqref{eq:double vertex derivative} implies $\delta L(z) \bar L^F(w) \sim - L^F(z) \, \overline{\delta L}(w)$, and thus
\begin{equation*}
L(z)\bar L(w) = L^F(z)\bar L^F(w) + L^F(z) \, \overline{\delta L}(w) +  \delta L(z)\bar L^F(w) + \delta L(z)\overline{\delta L}(w) \sim 0,
\end{equation*}
which means that the two Virasoro actions commute. 

It is easy to see that the formula \eqref{eq:Virasoro action 1} can be written as
\begin{equation*}
L(z) =  \mathcal Y\( \frac{(a_{-1})^2}{4\varkappa} \1_0 + \frac{\varkappa-1}{2\varkappa} a_{-2} \1_0 + \frac 1{\varkappa} \,\1_{-2},z\),
\end{equation*}
and similarly for \eqref{eq:Virasoro action 2}, which means that the vertex algebra structure is compatible 
with $\virCC$-module structure on $\tilde {\mathbb F}_\varkappa$.

We introduce the VOA structure in $\tilde{\mathbb F}_\varkappa$ by setting $\mathcal L(z) = L(z) + \bar L(z) =  L^F(z) + \bar L^F(z)$. 
One immediately checks that $\mathcal L(z)$ is a Virasoro quantum field with central charge $26,$ and satisfies $\mathcal L_{-1} \1_0 = 0.$ 
It suffices to check the remaining relation
\begin{equation}\label{eq:Virasoro L_{-1}}
[\mathcal L_{-1},\mathcal Y(v,z)] = \frac d{dz} \mathcal Y(v,z), \quad v \in \tilde{\mathbb F}_\varkappa,
\end{equation}
for each of the generating quantum fields, which is done by direct computations.
\end{proof}

\subsection{$\virCC$-module structure of $\tilde{\mathbb F}_\varkappa$ for generic $\varkappa$.}
We now describe the socle filtration of the $\virCC$-module $\tilde {\mathbb F}_\varkappa$ for generic $\varkappa$,
when it is completely analogous to the finite-dimensional and affine cases, given by Theorem \ref{thm:classical bimodule structure} and 
Theorem \ref{thm:affine bimodule structure}.
In this subsection we assume that 
$$\varkappa \notin \mathbb Q, \qquad c = 13 - \frac 6\varkappa - 6\varkappa,\qquad \bar c = 13 + \frac 6\varkappa + 6\varkappa.$$

For a $\vir_c$-module $\tilde V$, the restricted dual space $\tilde V'$ can be equipped with a $\vir_c$-action by
$$\<L_n \, v', v\> = \< v', L_{-n} v\>.$$
We denote the resulting dual module by $\tilde V^\star$.

We denote by $\tilde V_{\Delta,c}$ the irreducible $\vir_c$-module, generated
by a highest weight vector $\tilde v$ satisfying $L_0 \, \tilde v = \Delta \, \tilde v$ and $L_n \tilde v = 0$ for $n>0$.
For any $\lambda \in \h^*$, set
$$\Delta(\lambda) = \frac {\lambda(\lambda+2)}{4\varkappa} - \frac \lambda2 , \qquad 
\bar \Delta(\lambda) = -\frac {\lambda(\lambda+2)}{4\varkappa} - \frac \lambda2.$$

\medskip

\begin{thm}\label{thm:Virasoro bimodule structure}
There exists a filtration
\begin{equation}\label{eq:Virasoro filtration}
0 \subset \tilde{\mathbb F}_\varkappa^{(0)} \subset \tilde{\mathbb F}_\varkappa^{(1)} \subset \tilde{\mathbb F}_\varkappa^{(2)} = \tilde{\mathbb F}_\varkappa
\end{equation}
of $\virCC$-submodules of $\tilde{\mathbb F}_\varkappa$ such that
\begin{align}
\tilde{\mathbb F}_\varkappa^{(2)}/\tilde{\mathbb F}_\varkappa^{(1)} &\cong 
\bigoplus_{\lambda\in\mathbf P^+} \tilde V_{\Delta(-\lambda-2),c} \o \tilde V^\star_{\bar\Delta(-\lambda-2),\bar c} \label{eq:Virasoro socle F2},\\
\tilde{\mathbb F}_\varkappa^{(1)}/\tilde{\mathbb F}_\varkappa^{(0)} &\cong 
\bigoplus_{\lambda\in\mathbf P^+} \( \tilde V_{\Delta(\lambda),c} \o \tilde V^\star_{\bar\Delta(-\lambda-2),\bar c} \oplus 
\tilde V_{\Delta(-\lambda-2),c} \o \tilde V^\star_{\bar \Delta(\lambda),\bar c} \) \label{eq:Virasoro socle F1},\\
\tilde{\mathbb F}_\varkappa^{(0)} &\cong 
\bigoplus_{\lambda\in \mathbf P} \tilde V_{\Delta(\lambda),c} \o \tilde V^\star_{\bar\Delta(\lambda),\bar c} \label{eq:Virasoro socle F0}.
\end{align}
\end{thm}

\begin{proof}
One can derive from Proposition \ref{thm:Wakimoto to FF} that the for generic $\varkappa$ the reduction sends exact sequences
of $\ghat_k$-modules to exact sequences of $\vir_c$-modules, which implies in particular that
$$H_{DS}^n(\hat V_{\lambda,k}) = \begin{cases} \tilde V_{\Delta(\lambda),c}, & n=0 \\ 0, & n\ne 0 \end{cases},$$
It is then easy to check that the images $\tilde{\mathbb F}_\varkappa^{(0,1,2)}$ of the $\ghatKK$-submodules $\hat{\mathbb F}_\varkappa^{(0,1,2)}$ 
from Theorem \ref{thm:affine bimodule structure} under the two-sided quantum Drinfeld-Sokolov reduction satisfy the required properties.
\end{proof}

An alternative direct approach repeats the steps in the proof of Theorem \ref{thm:classical bimodule structure}.
In particular, we get a decomposition into blocks,
\begin{equation}
\tilde{\mathbb F}_\varkappa = \tilde{\mathbb F}_\varkappa(-1) \oplus \bigoplus_{\lambda\in \mathbf P^+} \tilde{\mathbb F}_\varkappa(\lambda).
\end{equation}

We also have the following Virasoro analogue of Corollary \ref{thm:classical positive subalgebra} and Theorem \ref{thm:affine positive subalgebra}.

\begin{thm}
There exists  a subspace $\tilde{\mathbf F}_\varkappa \subset \tilde{\mathbb F}_\varkappa$, satisfying
\begin{enumerate}
\item $\tilde{\mathbf F}_\varkappa$ is a vertex operator subalgebra of $\tilde{\mathbb F}_\varkappa$, and is 
generated by the quantum fields \eqref{eq:affine boson left}, \eqref{eq:affine boson right} and $\mathbb Y(1,z)$. 
In particular, $\tilde{\mathbf F}_\varkappa$ is a $\virCC$-submodule of $\tilde{\mathbb F}_\varkappa$.
\item As a $\virCC$-module, $\tilde{\mathbf F}_\varkappa$ is generated by the vectors $\{\1_\lambda\}_{\lambda \in \mathbf P^+}$, and we have
\begin{equation}
\label{eq:Virasoro positive decomposition}
\tilde{\mathbf F}_\varkappa \cong \bigoplus_{\lambda \in \mathbf P^+} \tilde V_{\Delta(\lambda),c} \o \tilde V^\star_{\bar\Delta(\lambda), \bar c}.
\end{equation}
\end{enumerate}
\end{thm}

\begin{proof}
The desired subspace $\tilde{\mathbf F}_\varkappa$ is the image of the vertex subalgebra 
$\hat{\mathbf F}_\varkappa$ under the two-sided quantum Drinfeld-Sokolov reduction. We leave technical details to the reader.
\end{proof}

\subsection{Semi-infinite cohomology of $\vir$}

The central charge for the diagonal action of $\vir$ in the modified regular representations  is equal to the special value 26.  
In this section we study the semi-infinite cohomology of $\vir$ with coefficients in $\tilde {\mathbb F}_\varkappa$ and in $\tilde {\mathbf F}_\varkappa$.

The properties of the appropriate "space of semi-infinite forms" $\boldsymbol{\tilde\Lambda}^\semiinfty$ for the Virasoro algebra 
are summarized in the following

\begin{prop}
Set $\boldsymbol{\tilde\Lambda}^\semiinfty = \bigwedge \vir_- \o \bigwedge \vir'_+$, where 
$\vir_- = \bigoplus_{n\le-2} \C L_n$ and $\vir_+ = \bigoplus_{n\ge-1} \C L_n$. Then

\begin{enumerate}
\item
The Clifford algebra, generated by $\{b_n, c_n \}_{n \in \Z}$ with relations
\begin{equation}
\label{eq:bc system relations}
\{b_m, b_n \} = \{ c_m, c_n \} = 0, \qquad
\{b_m, c_n\} = \delta_{m+n,0}.
\end{equation}
acts irreducibly on $\boldsymbol{\tilde\Lambda}^\semiinfty$, so that for any $\omega_- \in \bigwedge \vir_-, \, \omega_+ \in\bigwedge \vir'_+$ we have
\begin{align*}
b_n (1 \o \omega) &= 
\begin{cases} 0,& n \ge -1 \\ 1 \o (L_n\wedge\omega), & n\le-2 \end{cases}, \qquad
c_n (\omega \o 1) &= 
\begin{cases}  (L'_{-n}\wedge\omega) \o 1, & n\le 1 \\ 0,& n\ge 2 \end{cases}.
\end{align*}

\item $\boldsymbol{\tilde\Lambda}^{\semiinfty}$ is a bi-graded vertex superalgebra, with vacuum $\1 = 1 \o 1$,
generated by
\begin{alignat*}{3}
b(z) &= \sum_{n\in\Z} b_n z^{-n-2}, \qquad | b(z) | = -1, \quad  &\deg b(z) &= 2,\\
c(z) &= \sum_{n\in\Z} c_n z^{-n+1}, \qquad | c(z) | =  1, \quad  &\deg c(z) &= -1.
\end{alignat*}

\item $\boldsymbol{\tilde\Lambda}^\semiinfty$ has a $\vir$-module structure with central charge $c = -26$, defined by
$$\pi(L_n) = \sum_{m\in\Z} (m-n) :c_{-m} b_{n+m}:.$$
\end{enumerate}
\end{prop}

\begin{defn}
The BRST complex, associated with a $\vir$-module $\tilde V$ with central charge $c = 26$, is the complex
$C^{\semiinfty+\bullet} (\vir, \C\mathbf c;\tilde V) = \boldsymbol{\tilde\Lambda}^{\semiinfty+\bullet} \otimes \tilde V$, with the differential
\begin{equation}\label{eq:Virasoro differential}
\tilde{\mathbf d} = \sum_{n\in\Z} c_{-n} \pi_{\tilde V}(L_n) - \frac 12 \sum_{m,n\in \Z} (m-n) :c_{-m} c_{-n} b_{m+n}:.
\end{equation}
The corresponding cohomology is denoted $H^{\semiinfty+\bullet}(\vir,\C \mathbf c;\tilde V)$.
\end{defn}

As in the affine case, the special value $c = 26$ of the central charge is required to ensure that $\tilde{\mathbf d}^2 = 0$.

\begin{thm}
The vertex superalgebras $H^{\semiinfty+\bullet}(\vir,\C \mathbf c;\tilde{\mathbf F}_\varkappa)$ and $H^{\semiinfty+\bullet}(\vir,\C \mathbf c;\tilde{\mathbb F}_\varkappa)$ degenerate into the commutative superalgebras, and we have commutative algebra isomorphisms
\begin{equation}
\label{eq:Virasoro semi-infinite isomorphisms}
\begin{split}
H^{\semiinfty+\bullet}(\vir,\C \mathbf c;\tilde {\mathbf F}_\varkappa) \cong H^{\semiinfty+\bullet}(\ghat, \C\mathbf k;\hat{\mathbf F}_\varkappa), \qquad
H^{\semiinfty+\bullet}(\vir,\C \mathbf c;\tilde {\mathbb F}_\varkappa) \cong H^{\semiinfty+\bullet}(\ghat, \C\mathbf k;\hat{\mathbb F}_\varkappa).
\end{split}
\end{equation}
In particular, 
$$H^{\semiinfty+0}(\vir,\C \mathbf c;\tilde {\mathbf F}_\varkappa) \cong 
H^{\semiinfty+0}(\vir,\C \mathbf c;\tilde {\mathbb F}_\varkappa) \cong \C[\mathbf P]^W.$$
\end{thm}

\begin{proof}
The problem of computing the semi-infinite cohomology of $\vir$,
as well as its inherited algebra structure, has been extensively studied by mathematicians and physicists working in the string theory.
We take advantage of these results, and construct our proof by combining entire blocks from previous papers.

We note that for both $\tilde{\mathbf F}_\varkappa$ and $\tilde{\mathbb F}_\varkappa$ the diagonal action of $\vir$ does not contain 
additional vertex operator shifts, and is equal to the sum of two standard Feigin-Fuks actions.

The comprehensive answer for the cohomology of tensor products of Feigin-Fuks and/or irreducible modules was given in \cite{LZ2} 
for the most difficult case of the central charge $c = c_{p,q}$, corresponding to $\varkappa = \frac pq \in \mathbb Q$. 
Simplified (for the case of generic $\varkappa$) version of their computations, and the spectral sequence associated with filtrations 
of Theorem \ref{thm:Virasoro bimodule structure}, yield

$$H^{\semiinfty+n}(\vir,\C\mathbf c; \tilde{\mathbf F}_\varkappa(\lambda)) = 
\begin{cases}
\C, & n=0,3\\
0, & \text{otherwise}
\end{cases},$$
$$H^{\semiinfty+n}(\vir,\C\mathbf c; \tilde{\mathbb F}_\varkappa(-1)) = 
\begin{cases}
\C, & n=1,2\\
0, & \text{otherwise}
\end{cases},\qquad
H^{\semiinfty+n}(\vir,\C\mathbf c; \tilde{\mathbb F}_\varkappa(\lambda)) = 
\begin{cases}
\C, & n=0,2\\
\C^2, & n=1\\
0, & \text{otherwise}
\end{cases},$$
for each $\lambda \in \mathbf P^+$, as well as natural isomorphisms
$$H^{\semiinfty+0}(\vir,\C\mathbf c; \tilde{\mathbf F}_\varkappa(\lambda)) \cong H^{\semiinfty+0}(\vir,\C\mathbf c; \tilde{\mathbb F}_\varkappa(\lambda)).$$
The algebra structure of $H^{\semiinfty+0}(\vir,\C\mathbf c; \tilde{\mathbb F}_\varkappa)$ is in fact independent of $\varkappa$, as can be seen from the change of variables
$$p_n = \frac{a_n + \bar a_n}2, \qquad q_n = \frac {a_n - \bar a_n}{2\varkappa}.$$
Indeed, the new commutation relations become $[p_m,p_n] = [q_m,q_n] = 0$ and $[p_m,q_n] = \delta_{m+n,0}$,
and the diagonal Virasoro action is given by
$$\mathcal L(z) = :p(z)q(z): + p(z) - q(z).$$
For the special case $\varkappa = 1$, corresponding to the pairing of $c=1$ and $\bar c = 25$ modules, the cohomology of a bigger
vertex algebra $\mathcal A_{2D} = \bigoplus_{\lambda,\mu \in \Z} \tilde F_{\lambda,1} \o \tilde F_{\mu,-1}$ 
was identified in \cite{WZ} with the polynomial algebra $\C[x,y]$ in two variables. 
The subalgebra $H^{\semiinfty+0}(\vir,\C c; \tilde{\mathbb F}_\varkappa)$ is therefore isomorphic to the polynomial algebra $\C[\chi]$,
and we can take any nonzero cohomology class $\chi \in H^{\semiinfty+0}(\vir,\C c; \tilde{\mathbb F}_\varkappa(1))$ as the generator.

The vertex superalgebra structures on $H^{\semiinfty+\bullet}(\vir,\C\mathbf c; \tilde{\mathbf F}_\varkappa)$ and $H^{\semiinfty+\bullet}(\vir,\C\mathbf c; \tilde{\mathbb F}_\varkappa)$ degenerate into commutative superalgebras. It is clear that both are free $\C[\chi]$-modules.

It follows immediately that $H^{\semiinfty+\bullet}(\vir,\C\mathbf c; \tilde{\mathbf F}_\varkappa) \cong \C[\chi] \otimes \bigwedge^\bullet[\eta]$,
where we can pick any non-zero element $\eta \in H^{\semiinfty+3}(\vir,\C\mathbf c; \tilde{\mathbf F}_\varkappa(0))$. This settles the
case of $\tilde{\mathbf F}_\varkappa$.

To get the generators of $H^{\semiinfty+\bullet}(\vir,\C\mathbf c; \tilde{\mathbb F}_\varkappa)$, we pick non-zero representatives
$$\xi_{-1} \in H^{\semiinfty+1}(\vir,\C\mathbf c; \tilde{\mathbb F}_\varkappa(-1)), \qquad
\eta_0 \in H^{\semiinfty+1}(\vir,\C\mathbf c; \tilde{\mathbb F}_\varkappa(0)),$$
such that $\eta_0$ is not proportional to $\chi \cdot \xi_{-1}$. 
One can check that $\eta_0 \xi_{-1} \ne 0$, and as in Theorem \ref{thm:classical cohomology}
it follows that $H^{\semiinfty+\bullet}(\vir,\C\mathbf c; \tilde{\mathbb F}_\varkappa) \cong \C[\chi] \otimes \bigwedge^\bullet[\xi_{-1},\eta_0]$.
The statement now follows from Theorem \ref{thm:classical cohomology} and Corollary \ref{thm:affine cohomology}.
\end{proof}

\begin{rem}
It would be nice to get a direct proof of isomorphisms \eqref{eq:Virasoro semi-infinite isomorphisms}
by using the techniques of the quantum Drinfeld-Sokolov reduction.
\end{rem}

\section{Extensions, generalizations, conjectures}

\subsection{Heterogeneous vertex operator algebra}
As we mentioned above, the vertex algebra construction for the Virasoro algebra can be obtained from their affine analogues by
applying the two-sided quantum Drinfeld-Sokolov reduction to the left and right $\ghat$-action. One can consider a similar construction 
where the reduction is only applied to the affine action on one side, thus leading to a vertex operator algebra with
two commuting actions of  $\ghat_k$ and $\vir_{\bar c}$ with appropriate $k,\bar c$. Indeed, one can see that the following
gives a direct realization of such a vertex algebra.

\begin{thm}\label{thm:dual pair action}
Let $\varkappa \ne 0$. Set $k = \varkappa - h^\vee,\  \bar c = 13 + 6\varkappa + \frac 6\varkappa$, and $\check {\mathbb F}_\varkappa = \hat F(\beta,\gamma) \o \tilde{\mathbb F}_\varkappa.$ Then
\begin{enumerate}
\item
The space $\check{\mathbb F}_\varkappa$ has a $\ghat_k \oplus \vir_{\bar c}$-module structure, defined by
\begin{equation}\label{eq:heterogeneous affine current}
\begin{split} 
\mathbf e(z) &= \gamma(z),\\
\mathbf h(z) &= 2:\beta(z) \gamma(z): + a(z),\\
\mathbf f(z) &= -:\beta(z)^2 \gamma(z): - \beta(z) a(z) - k \beta'(z) - \mathbb Y(-2,z),
\end{split}
\end{equation}
\begin{equation} \label{eq:heterogeneous Virasoro current}
\bar L(z) = -\frac {:\bar a(z)^2:}{4\varkappa}  + \frac{\varkappa+1}{2\varkappa} \bar a'(z)
 + \frac 1\varkappa \, \mathbb Y(-2,z) \gamma(z).
\end{equation}
\item
The space $\check{\mathbb F}_\varkappa$ has a compatible VOA structure with $\rank \check{\mathbb F}_\varkappa = 28$.
\end{enumerate}
\end{thm}
\begin{proof}
The verification of commutation relations is straightforward. 
We define the Virasoro quantum field by
\begin{equation}
\label{eq:heterogeneous total Virasoro}
\mathcal L(z) = \frac 1{2\varkappa} \(\frac {:\mathbf h(z)^2:}2 + :\mathbf e(z)\mathbf f(z): + :\mathbf f(z)\mathbf e(z):\)
 + \frac {\mathbf h'(z)}2 + \bar L(z)
\end{equation}
The central charge for the Sugawara construction modified by $\frac {\mathbf h'(z)}2$ is equal to $\frac {3k}{k+h^\vee}-6k$, and we compute
$$\rank \check{\mathbb F}_\varkappa = \( \frac{3(\varkappa-h^\vee)}\varkappa - 6(\varkappa-h^\vee) \) + \(13 + \frac 6\varkappa + 6\varkappa \) = 28.$$
\end{proof}

We will call the vertex operator algebra of Theorem \ref{thm:dual pair action} the heterogeneous VOA.
Note (see \cite{Li} and references therein) that the central charge $c=28$ appears as the critical value in the study of 2D gravity in the light-cone gauge!


The structure of the bimodule $\check {\mathbb F}_\varkappa$ in the generic case is again quite similar to the non-semisimple bimodule $\mathfrak R(G_0)$.
From now on we assume that

$$\varkappa \notin \mathbb Q, \qquad k = \varkappa - h^\vee, \qquad \bar c = 13 + \frac 6\varkappa + 6\varkappa.$$

\begin{thm}\label{thm:dual pair structure}
There exists a filtration
$$0 \subset \check{\mathbb F}_\varkappa^{(0)} \subset \check{\mathbb F}_\varkappa^{(1)} \subset \check{\mathbb F}_\varkappa^{(2)} = 
\check{\mathbb F}_\varkappa$$
of $\ghat_k \oplus \vir_{\bar c}$-submodules of $\check{\mathbb F}_\varkappa$, such that
\begin{align}
\check{\mathbb F}_\varkappa^{(2)}/\check{\mathbb F}_\varkappa^{(1)} &\cong 
\bigoplus_{\lambda\in\mathbf P^+} \hat V_{-\lambda-2,k} \o \tilde V^\star_{\bar\Delta(-\lambda-2),\bar c}, \\
\check{\mathbb F}_\varkappa^{(1)}/\check{\mathbb F}_\varkappa^{(0)} &\cong 
\bigoplus_{\lambda\in\mathbf P^+} \( \hat V_{\lambda,k} \o \tilde V^\star_{\bar\Delta(-\lambda-2),\bar c} \oplus 
\hat V_{-\lambda-2,k} \o \tilde V^\star_{\bar\Delta(\lambda),\bar c} \), \\
\check{\mathbb F}_\varkappa^{(0)} &\cong 
\bigoplus_{\lambda\in\mathbf P} \hat V_{\lambda,k} \o \tilde V^\star_{\bar\Delta(\lambda),\bar c}.
\end{align}
\end{thm}

The heterogeneous VOA contains a vertex operator subalgebra, analogous to the classical Peter-Weyl subalgebra $\mathfrak R(G) \subset \mathfrak R(G_0)$.

\begin{thm}
There exists  a subspace $\check{\mathbf F}_\varkappa \subset \check{\mathbb F}_\varkappa$, satisfying
\begin{enumerate}
\item $\check{\mathbf F}_\varkappa$ is a vertex operator subalgebra of $\check{\mathbb F}_\varkappa$, and is generated by the fields \eqref{eq:heterogeneous affine current}, \eqref{eq:heterogeneous Virasoro current}, and $\mathbb Y(1,z)$. In particular, $\check{\mathbf F}_\varkappa$ is a $\ghat_k \oplus \vir_{\bar c}$-submodule of $\check{\mathbb F}_\varkappa$.
\item As a $\ghat_k \oplus \vir_{\bar c}$-module, $\check{\mathbf F}_\varkappa$ is generated by the vectors $\{\1_\lambda\}_{\lambda \in \mathbf P^+}$, and we have
\begin{equation}
\check{\mathbf F}_\varkappa \cong \bigoplus_{\lambda \in \mathbf P^+} \hat V_{\lambda,k} \o \tilde V_{\bar\Delta(\lambda),\bar c}.
\end{equation}
\end{enumerate}
\end{thm}

The proofs of the above theorems are obtained from their affine counterparts by applying the quantum Drinfeld-Sokolov
reduction to (right) $\ghat_{\bar k}$-action, similarly to the Virasoro case.

The fact that the ranks of VOAs $\check{\mathbb F}_\varkappa$ and $\check{\mathbf F}_\varkappa$ are equal to 28 naturally
leads to the consideration of the semi-infinite cohomology of these modules. We note that although the total Virasoro
quantum field $\mathcal L(z)$ does not commute with $\ghat$, the spaces $\check{\mathbb F}_\varkappa$ and 
$\check{\mathbf F}_\varkappa$ can be regarded as modules over the semi-direct product $\virghat$, such that
\begin{equation}
\begin{split}
[\mathcal L_m, \mathbf e_n] &= -(n+m+1)\, \mathbf e_{m+n},\\
[\mathcal L_m, \mathbf f_n] &= (m-n+1)\, \mathbf f_{m+n},\\
[\mathcal L_m, \mathbf h_n] &= -n\, \mathbf h_{m+n} + m(m+1)\delta_{m+n,0}\, k.
\end{split}
\end{equation}

The semi-infinite cohomology is defined by the BRST complex of the subalgebra $\virnil$, where as before
$\nhat_+ = \bigoplus_{n \in \Z} \C \mathbf e_n$ is the nilpotent loop subalgebra of $\ghat$. The condition $c=28$
ensures that the differential squares to zero.

\begin{defn} The BRST complex, associated with a $\virnil$-module $\check V$ with Virasoro central charge $c=28$ is the complex
$C^{\semiinfty+\bullet} (\vir, \C\mathbf c;\check V) = \boldsymbol{\tilde\Lambda}^{\semiinfty+\bullet} \o \hat\Lambda(\psi,\psi^*) \o \check V$, with the differential
\begin{equation}\label{eq:heterogeneous differential}
\begin{split}
\check{\mathbf d} &= \sum_{n\in\Z} c_{-n} \pi_{\check V}(\mathcal L_n) - \frac 12 \sum_{m,n\in \Z} (m-n) :c_{-m} c_{-n} b_{m+n}: + \\
&+ \sum_{n \in \Z} \psi^*_{-n} \pi_{\check V}(\mathbf e_n) - \sum_{m,n\in\Z} -(m+n+1) c_{-m} :\psi^*_{-n} \psi_{m+n}:
\end{split}
\end{equation}
The corresponding cohomology is denoted $H^{\semiinfty+\bullet}(\virnil,\C \mathbf c;\check V)$.
\end{defn}

\begin{prop} The vertex algebras $H^{\semiinfty+\bullet}(\virnil,\C\mathbf c; \check{\mathbf F}_\varkappa)$ and
$H^{\semiinfty+\bullet}(\virnil,\C\mathbb c; \check{\mathbf F}_\varkappa)$ degenerate into commutative algebra structures, and
we have isomorphisms
$$H^{\semiinfty+\bullet}(\virnil,\C\mathbf c; \check{\mathbf F}_\varkappa) \cong H^{\semiinfty+\bullet}(\vir,\C\mathbf c; \tilde{\mathbf F}_\varkappa),$$
$$H^{\semiinfty+\bullet}(\virnil,\C\mathbf c; \check{\mathbb F}_\varkappa) \cong H^{\semiinfty+\bullet}(\vir,\C\mathbf c; \tilde{\mathbb F}_\varkappa).$$
In particular,
$$H^{\semiinfty+0}(\virnil,\C\mathbf c; \check{\mathbf F}_\varkappa) \cong 
H^{\semiinfty+0}(\virnil,\C\mathbf c; \check{\mathbb F}_\varkappa) \cong \C[\mathbf P]^W.$$
\end{prop}

\begin{proof}
We use the technique from \cite{Li}, where similar isomorphisms were established for relative cohomology spaces.
Let $\mathbf d_{\n_+} = \sum_{n \in \Z} \psi^*_{-n} \pi(\mathbf e_n)$ be the BRST differential for $\n_+$; 
one can show that $\mathbf d_{\n_+}^2 = 0 = \{\mathbf d_{\n_+},\check{\mathbf d}\}$, which leads to the spectral 
sequence associated with decomposition $\check{\mathbf d} = \mathbf d_{\n_+} + (\check{\mathbf d} - \mathbf d_{\n_+})$.
Computing the cohomology with respect to $\mathbf d_{\n_+}$ first, and using Proposition \ref{thm:Wakimoto to FF},
we get the desired statement. For full technical details (there is a slight difference between BRST reduction
for $\n_+$ and quantum Drinfeld-Sokolov reduction, but it doesn't affect the outcome) we refer the reader to \cite{Li}.
\end{proof}

\subsection{General construction of vertex operator algebras and equivalence of categories}

The vertex operator algebras constructed in the previous sections can be built, like conformal field theories, by pairing the left and right modules 
from certain equivalent categories of representations of infinite-dimensional Lie algebras. The operators $\mathcal Y(\cdot,z)$ are then constructed 
by pairing the left and right intertwining operators. There is a unique choice of the structural coefficients for such pairing that would ensure 
the locality condition for the vertex operator algebras; these coefficients are determined by the tensor structure on the category of representations. 
Conversely, a natural VOA structure on a bimodule can be used to establish the equivalence of the left and right tensor categories.

The vertex operator algebra constructions in this paper deal with the pairings of different categories of modules. In the affine case, 
we pair the $\ghat$ modules on levels $k$ and $\bar k = -2h^\vee - k$, symmetric with respect to the critical level $-h^\vee$;
the equivalence of the corresponding tensor categories was studied in \cite{Fi}. In the Virasoro case,
we pair the modules with central charges $c$ and $\bar c = 26-c$.

The theorems of Peter-Weyl type can be extended to the quantum group $\mathcal U_q(\g)$, associated with $G$. 
On one hand, the modules from the category $\mathcal O$ can be $q$-deformed into modules over $\mathcal U_q(\g)$;
on the other hand, one can define $q$-deformations $\mathfrak R_q(G)$ and $\mathfrak R_q(G_0)$ of the algebras of regular functions, \
which have especially simple description for $\g = \sl(2,\C)$. When $q$ is not a root of unity, 
we have the quantum analogues of isomorphisms \eqref{eq:Peter-Weyl classical}, \eqref{eq:Peter-Weyl projective}.

The Drinfeld-Kohno theorem establishes an isomorphism of tensor categories of representations of the quantum groups and affine Lie algebras when 
$q = \exp\(\frac{\pi i}{k+h^\vee}\)$. This equivalence, also extended to $\mathcal W$-algebras, was made explicit in \cite{S},
where intertwining operators for $\mathcal U_q(\g)$ were directly identified with their VOA counterparts for $\ghat_k$ and $\mathcal W(\ghat_k)$;
the key ingredients were the geometric results in \cite{V} on the homology of configuration spaces.
The construction in \cite{S} built conformal field theories, associated to affine Lie algebras 
and $\mathcal W$-algebras based on their quantum group counterparts, and can be modified to produce the vertex algebras discussed in this paper.

The Drinfeld-Kohno equivalence also allows to couple categories of different types, producing in particular the heterogeneous VOA of the 
previous subsection. Another important case is the Frenkel-Kac construction, which corresponds to the pairing of modules for $\g$ and $\mathcal W(\ghat)$
with central charge $c=\dim \h$ (see \cite{F}). However, in general pairings between the quantum group and the affine Lie (or $\mathcal W$-) algebra
lead to the generalized vertex algebra structures, satisfying a braided version of the commutativity axiom. An example of such structure
was proposed in \cite{MR}.

\subsection{Integral central charge, semi-infinite cohomology and Verlinde algebras}

In this work we studied the structure of the generalized Peter-Weyl bimodules for $\ghat$ only for the generic values $k \notin \mathbb Q$ 
of the central charge. The structure of these bimodules when $k$ is integral is more complex and undoubtedly even more interesting. 
In the most special case when $k = \bar k = - h^\vee$, we get a regular representation of the affine Lie algebra $\ghat$ at the critical level, 
which can be viewed as the direct counterpart of the finite-dimensional case. This space admits a realization as a certain space of 
meromorphic functions on the affine Lie group $\hat G$, and subspaces of spherical functions with respect to conjugation give rise to solutions 
of the quantum elliptic Calogero-Sutherland system, generalizing the trigonometric analogue in the finite-dimensional case \cite{EFK}.

Another special case is $k = - h^\vee+1, \bar k = -h^\vee - 1$, when the quantum group degenerates into its classical counterpart.
In this case the left and right Fock spaces used in our construction each have separate vertex algebra structures, and the 
operators $\mathbb Y(-2,w)$, which play an important role in this paper, are factored into products of left and right vertex operators used 
in the basic representations of $\ghat$. The semi-infinite cohomology of the corresponding $\mathcal W$-algebras is fundamental to the 
string theory, and was studied in \cite{WZ} for the Virasoro algebra, and in \cite{BMP} for $\mathcal W_3$.

For positive integral $k$ one expects the existence of truncated versions $\hat{\mathbf F}_{k+h^\vee}$ and $\hat{\mathbb F}_{k+h^\vee}$ of our vertex operators algebras, 
similar to the truncation in the conformal field theory, where the positive dominant cone $\mathbf P^+$ is replaced by the alcove 
$\mathbf P^+_k \subset \mathbf P^+$. Then the relative semi-infinite cohomology of $\ghat$ with coefficients in $\hat{\mathbf F}_{k+h^\vee}$ and 
$\hat{\mathbb F}_{k+h^\vee}$ with respect to the center should be truncated correspondingly. 
The identification of the zero semi-infinite cohomology groups with the representation ring of $G$ in Corollary \ref{thm:affine cohomology} leads to the following conjecture.

\begin{conj}
\label{thm:conjecture} 
For positive integral $k$, let $\mathbf V_k(\ghat)$ denote the Verlinde algebra associated with integrable level $k$ representations of $\ghat$, and let 
$\mathbb V_k(\ghat)$ denote its counterpart associated to the big projective modules (see \cite{La}). Then 
we have commutative algebra isomorphisms
$$H^{\semiinfty+0}(\ghat,\C\mathbf k; \hat{\mathbf F}_{k+h^\vee})  \cong \mathbf V_k(\ghat), \qquad \qquad
H^{\semiinfty+0}(\ghat,\C\mathbf k; \hat{\mathbb F}_{k+h^\vee})  \cong \mathbb V_k(\ghat).$$
\end{conj}

In other words, the most essential part of the VOA structure, embodied in the 0th cohomology, is equivalent to the structure of the fusion rules 
of the tensor category of $\ghat$-modules, encoded in the Verlinde algebra.

It was also realized recently (see \cite{FHT} and references therein) that the Verlinde algebra $\mathbf V_k(\ghat)$ admits an
 alternative realization in terms of twisted equivariant K-theory ${}^{k+h^\vee}K_G^{\dim G}(G)$ of a compact simple Lie group $G$
 (which, in the notations of \cite{FHT}, is the compact form of the complex Lie group which we denoted by $G$ in this paper).

Thanks to the results of \cite{FHT}, the first isomorphism of our conjecture can be restated in a more invariant form.
\begin{conj}
We have a natural commutative algebra isomorphism
$$H^{\semiinfty+0}(\ghat,\C\mathbf k; \mathfrak R'_{k+h^\vee}(\hat G)) \cong {}^{k+h^\vee}K_G^{\dim G}(G).$$
\end{conj}
A conceivable direct geometric proof of the last isomorphism might combine the realization of the left 
hand side using the works \cite{GMS} and \cite{AG} with the interpretation of the right hand side given in the works \cite{FHT} and \cite{AS}.
A similar K-theoretic interpretation of $H^{\semiinfty+0}(\ghat,\C\mathbf k; \hat{\mathbb F}_{k+h^\vee})$ 
in our second conjecture might add another twirl to the twisted equivariant K-theory.


\begin{thebibliography}{}
\frenchspacing

\bibitem[AG]{AG}
S. Arkhipov, D. Gaitsgory,
Differential operators on the loop group via chiral algebras. Int. Math. Res. Not. 2002, no. 4, 165-210

\bibitem[AS]{AS}
M. Atiyah, G. Segal.
Twisted K-theory, math.KT/0407054


\bibitem[BF]{BF}
D. Bernard, G. Felder, 
Fock representations and BRST cohomology in ${\rm SL}(2)$ current algebra. 
Comm. Math. Phys. 127 (1990), no. 1, 145--168.

\bibitem[BGG]{BGG}
J. Bernstein, I. Gelfand, S. Gelfand,
A certain category of $\g$-modules.
(Russian) Funkcional. Anal. i Prilo\v zen. 10 (1976), no. 2, 1--8. 

\bibitem[BMP]{BMP}
P. Bouwknegt, J. McCarthy, K. Pilch,
The ${\mathcal W}_3$ algebra. Modules, semi-infinite cohomology and BV algebras. 
Lecture Notes in Physics. Monographs, 42. Springer-Verlag, Berlin, 1996.

\bibitem[EFK]{EFK}
P. Etingof, I. Frenkel, A. Kirillov Jr.,
Spherical functions on affine Lie groups.
Duke Math. J. 80 (1995), no. 1, 59--90. 

\bibitem[Fe]{Fe}
B. Feigin,
Semi-infinite homology of Lie, Kac-Moody and Virasoro algebras.
(Russian) Uspekhi Mat. Nauk 39 (1984), no. 2(236), 195--196.

\bibitem[FeFr1]{FeFr}
B. Feigin, E. Frenkel, 
Representations of affine Kac-Moody algebras and bosonization.
Physics and mathematics of strings, 271--316, World Sci. Publishing, Teaneck, NJ, 1990.

\bibitem[FeFr2]{FeFrDS}
B. Feigin, E. Frenkel, 
Affine Kac-Moody algebras at the critical level and Gelfand-Dikii algebras.
Infinite analysis, Part A, B (Kyoto, 1991), 197--215,
Adv. Ser. Math. Phys., 16,
World Sci. Publishing, River Edge, NJ, 1992. 

\bibitem[FeFu]{FeFu}
B. Feigin, D. Fuks, 
Verma modules over a Virasoro algebra.
Funktsional. Anal. i Prilozhen. 17 (1983), no. 3, 91--92. 

\bibitem[FeP]{FP}
B. Feigin, S. Parkhomenko,
Regular representation of affine Kac-Moody algebras. Algebraic and geometric methods in mathematical physics (Kaciveli, 1993), 
415--424, Math. Phys. Stud., 19, Kluwer Acad. Publ., Dordrecht, 1996.

\bibitem[Fi]{Fi}
M. Finkelberg, An equivalence of fusion categories. Geom. Funct. Anal. 6 (1996), no. 2, 249--267.

\bibitem[FHT]{FHT}
D. Freed, M. Hopkins, C. Teleman,
Twisted K-theory and loop group representations I.
math.AT/0312155

\bibitem[FrB]{BFr}
E. Frenkel, D. Ben-Zvi,
Vertex algebras and algebraic curves. Mathematical Surveys and Monographs, 88. American Mathematical Society, Providence, RI, 2001

\bibitem[F]{F}
I. Frenkel,
Representations of Kac-Moody algebras and dual resonance models,
Applications of group theory in physics and mathematical physics (Chicago, 1982), 325--353, 

\bibitem[FGZ]{FGZ}
I. Frenkel, H. Garland, G. Zuckerman,
Semi-infinite cohomology and string theory. 
Proc. Nat. Acad. Sci. U.S.A. 83 (1986), no. 22, 8442--8446.

\bibitem[FLM]{FLM}
I. Frenkel, J. Lepowsky, A. Meurman
Vertex operator algebras and the Monster. Pure and Applied Mathematics, 134. Academic Press, Inc., Boston, MA, 1988.

\bibitem[FZ]{FZ}
I. Frenkel, Y. Zhu,
Vertex operator algebras associated to representations of affine and Virasoro algebras. 
Duke Math. J. 66 (1992), no. 1, 123--168.

\bibitem[GMS]{GMS}
V. Gorbounov, F. Malikov, V. Schechtman,
On chiral differential operators over homogeneous spaces.  Int. J. Math. Math. Sci.  26  (2001),  no. 2, 83--106.

\bibitem[HS]{HS}
P. Hilton, U. Stammbach,
A course in homological algebra. Graduate Texts in Mathematics, Vol. 4. Springer-Verlag, New York-Berlin, 1971.

\bibitem[La]{La}
A. Lachowska,
A counterpart of the Verlinde algebra for the small quantum group. Duke Math. J. 118 (2003), no. 1, 37--60. 

\bibitem[L]{Li}
B. Lian, Semi-infinite homology and 2D quantum gravity.
PhD thesis, Yale University, (1991).

\bibitem[LZ1]{LZ}
B. Lian, G. Zuckerman,
New perspectives on the BRST-algebraic structure of string theory.
Comm. Math. Phys. 154 (1993), no. 3, 613--646.

\bibitem[LZ2]{LZ2}
B. Lian, G. Zuckerman,
Semi-infinite homology and $2$D gravity. I.  
Comm. Math. Phys.  145  (1992),  no. 3, 561--593.

\bibitem[MR]{MR}
G. Moore, N. Reshetikhin, 
A comment on quantum group symmetry in conformal field theory. 
Nuclear Phys. B 328 (1989), no. 3, 557--574.

\bibitem[PS]{PS}
A. Pressley, G. Segal,
Loop groups. Oxford Mathematical Monographs. Oxford Science Publications. 
The Clarendon Press, Oxford University Press, New York, 1986.


\bibitem[S]{S}
K. Styrkas, Quantum groups, conformal field theories, and duality of tensor categories.
PhD thesis, Yale University, (1998).

\bibitem[V]{V}
A. Varchenko,
Multidimensional hypergeometric functions and representation theory of Lie algebras and quantum groups.
Advanced Series in Mathematical Physics, 21. World Scientific Publishing Co., Inc., River Edge, NJ, 1995

\bibitem[W]{W}
F. Williams,
The cohomology of semisimple Lie algebras with coefficients in a Verma module.
Trans. Amer. Math. Soc.  240  (1978), 115--127. 

\bibitem[WZ]{WZ}
E.Witten, B.Zwiebach,
Algebraic structures and differential geometry in 2d string theory.
Nucl. Phys. B377 (1992) , 55--112
\end{thebibliography}
\end{document}